\documentclass[smallextended,referee,envcountsect,]{svjour3}
\smartqed
\usepackage{geometry}
\usepackage{rotating}
\usepackage{graphicx}
\usepackage{float}
\usepackage{amsmath,amssymb,amsfonts,amstext,amsmath}
\usepackage{epsfig}
\usepackage{setspace}
\usepackage{hhline}
\usepackage{color,soul}
\usepackage[dvipsnames]{xcolor}
\definecolor{Light}{gray}{.90}
\usepackage[justification=centering]{caption}
\usepackage{multirow}
\usepackage{hyperref}
\usepackage{booktabs}
\usepackage{bm}
\usepackage[english]{babel}
\usepackage{tikz}
\usepackage{tikzscale}
\usepackage{epstopdf}
\usepackage{adjustbox}
\usepackage{longtable}
\usepackage{accents}
\usepackage{lineno}
\usepackage{multirow}
\usepackage{colortbl}
\usepackage[linesnumbered,ruled,vlined]{algorithm2e}
\usepackage{ragged2e}
\usepackage{amsmath,amssymb}
\makeatletter
\usepackage{breqn}
\usepackage{pdflscape}
\usepackage{afterpage}
\usepackage{capt-of}
\usepackage{dsfont}
\usepackage{gensymb}
\newtheorem{assumption}{Assumption}

\newcommand{\iter}[1]{\ensuremath{\langle #1 \rangle}}
\newcommand{\itrpt}[3]{\ensuremath{\bm{#1}^{\iter{#2}}_{#3}}}

\newcommand{\idxpt}[3]{\ensuremath{\bm{#1}^{(#2)}_{#3}}}
\newcommand{\idxptT}[3]{\ensuremath{\bm{#1}^{(#2)\top}_{#3}}}
\newcommand{\idxval}[3]{\ensuremath{#1^{(#2)}_{#3}}}

\newcommand{\tidxpt}[3]{\ensuremath{\tilde{\bm{#1}}^{(#2)}_{#3}}}

\journalname{JOTA}

\begin{document}

\title{SAGE: A Set-based Adaptive Gradient Estimator}


\author{Lorenzo Sabug Jr. \and Fredy Ruiz \and Lorenzo Fagiano}

\institute{Lorenzo Sabug, Jr. \at
           Dipartimento di Elettronica, Informazione e Bioingegneria, Politecnico di Milano, \\Piazza Leonardo da Vinci 32, 20133 Milano, Italy \\
           \textit{Present address:} Department of Electrical and Electronic Engineering, Imperial College, \\Exhibition Road, SW7 2AZ London, United Kingdom. \email{l.sabug21@imperial.ac.uk}
           \and
           Fredy Ruiz \at
           Dipartimento di Elettronica, Informazione e Bioingegneria, Politecnico di Milano, \\Piazza Leonardo da Vinci 32, 20133 Milano, Italy. \email{fredy.ruiz@polimi.it}
           \and
           Lorenzo Fagiano \at
           Dipartimento di Elettronica, Informazione e Bioingegneria, Politecnico di Milano, \\Piazza Leonardo da Vinci 32, 20133 Milano, Italy. \email{lorenzo.fagiano@polimi.it}
}

\date{Received: date / Accepted: date}

\maketitle

\begin{abstract}
We introduce a set membership perspective on gradient estimation for black-box scalar functions, in which the gradient is not approximated as a single point estimate but characterised as a member of a set of admissible gradients consistent with available function samples and regularity assumptions. Bounds on gradient uncertainty, obtained from standard multivariate Taylor series results, are used to represent this admissible set through linear inequalities, making gradient uncertainty explicit and directly exploitable. Building on this formulation, we propose an adaptive strategy to refine the admissible gradient set to a desired precision by selectively acquiring new samples only when necessary. The resulting framework naturally accommodates function evaluations affected by noise with finite but \textit{a priori} unknown bounds, yields the theoretical best attainable precision of the gradient set, and provides the optimal sampling distance from the point of interest to maximise information gain during refinement. Within this framework, we introduce the Set-Based Adaptive Gradient Estimator (SAGE), which extracts an estimate of the true gradient from the admissible set whilst preserving high sample efficiency and robustness to noise. Numerical experiments on large collections of random gradient-query points, as well as in the context of a first-order optimisation method, compared the performance of SAGE against commonly used gradient estimators from the literature. Our results show that SAGE achieved competitive gradient accuracy across noise levels without requiring any \textit{a priori} noise characterisation. In a first-order optimisation setting, noise-informed competitors can outperform SAGE when such characterisation is available. Nevertheless, SAGE remained broadly competitive, making it particularly attractive in settings where noise characterisation is impractical to obtain.
\end{abstract}


\keywords{Gradient estimation \and Data-driven methods \and Set membership approach}
\subclass{03E75 \and  90C56}


\section{Introduction}

The gradient of a scalar function is an essential ingredient in many optimisation \cite{papadimitriou2008aerodynamic} and estimation \cite{wan2000unscented} workflows, e.g., in machine learning, scientific computing, and multidisciplinary design. However, in many applications, like in experiment- or simulation-based design and/or optimisation, the closed-form expression of the function is unavailable, and individual function values can only be acquired by performing a resource-intensive (e.g., expensive or time-consuming) evaluation at discrete points in the search space, i.e., the function is a ``black-box''. As a result, information on the gradient is unavailable, or otherwise difficult to compute. As a workaround, gradients are estimated by acquiring additional samples close to the point of interest (which we refer to as the ``query point''), and computing the directional slopes. However, acquiring these additional samples may be similarly time-consuming or expensive (e.g., when real-world experiments are needed), and such finite-difference approaches can be susceptible to noisy samples. Hence, developing gradient estimation approaches which are both sample efficient and robust to noise is of interest and an active area of research. 

Gradient estimation using zeroth-order information (as for derivative-free optimisation, or DFO) is a mature topic, with well-developed theoretical foundations, practical approaches, and quantification of error bounds \cite{audet2017derivative,conn2009introduction}. Classical finite-difference-based methods, usually forward finite differences (FFD) and central finite differences (CFD), use additional (``auxiliary'') samples at coordinate-wise perturbations from the current query point to estimate directional derivatives, and are widely used in engineering design and scientific computing. In the DFO literature, model-based methods construct surrogate models and use their gradients as the estimates to the true one. For instance, simplex gradients and linear interpolation gradients fit an affine model to a poised set of samples, and the resulting gradient estimation error can be bounded by the sampling radius and geometric conditioning under Lipschitz gradient assumptions \cite{conn2009introduction,regis2014calculus}. The analysis of generalised simplex gradients and other variants have tackled accuracy for over-/underdetermined sampling and different samples set geometries \cite{hare2022error}.

Smoothing-based estimators have also emerged in significance, particularly for large-scale and noisy settings, because of their simplicity and convergence guarantees in expectation when deployed in DFO. These include Gaussian smoothing and two-point gradient estimators, approximating gradients using random auxiliary sample perturbations and function differences \cite{berahas2022comparison,duchi2015optimal,nesterov2017random}. Noise-robust gradient estimation has also been studied through auxiliary sample step size selection and error bounds for FD estimators under bounded evaluation noise \cite{more2012estimating,scheinberg2022finite}. In addition, practical variants improve noise robustness by collecting FD estimates across multiple step sizes \cite{boresta2022mixed}, and by reusing samples within local neighbourhoods \cite{wang2024relizo}. Linear interpolation was shown in \cite{berahas2019linear} to yield competitive gradients compared to some smoothing estimators under appropriate geometric conditions.

A unifying feature across most of the above-discussed approaches is the \textit{estimator-first} paradigm: a point-valued gradient estimate is defined \textit{a priori} (from FD, simplex/interpolation gradients, etc.), and its error is then analysed under smoothness assumptions, noise level, and sampling geometry \cite{audet2017derivative,berahas2022comparison,conn2009introduction}. In contrast, this work discusses a \textit{set-identification} viewpoint: rather than committing \textit{a priori} to a single point estimate, we derive inequality constraints on the unknown gradient, induced by existing samples and smoothness assumptions. From these, we construct a set of admissible gradients guaranteed to contain the true one. A point gradient is then selected from this set, and the geometry of the admissible gradients set is explicitly used to guide adaptive sampling only when necessary to achieve a desired gradient accuracy. Based on our results, we propose a new approach, named the Set-Based Adaptive Gradient Estimator (SAGE), which computes an estimate of the gradient inside the set of admissible ones, and intelligently acquires new samples only as needed to improve this estimate. SAGE features both high sample efficiency and high robustness to noise, as we also demonstrate in comparative statistical tests with other methods (code available in \url{https://github.com/lorenzosabugjr/SAGE}). Summarising, we put forth the following main contributions:

\begin{enumerate}
    \item We introduce a set membership paradigm for gradient estimation using zeroth-order information, in which the gradient is treated as an unknown variable constrained by sample-induced Taylor-based inequalities. Rather than constructing a point estimator, we identify a polyhedral admissible gradient region that is guaranteed to contain the true gradient under stated smoothness assumptions, enabling explicit representation of gradient uncertainty.
    
    \item In the practically relevant case of \textit{a priori} unknown Hessian norms and Lipschitz constants, we estimate both the admissible gradient set and the required smoothness parameters directly from data via linear programming. The geometry of the admissible set (e.g., diameter or unboundedness) is explicitly exploited to trigger adaptive sampling, ensuring that additional function evaluations are performed only when necessary.
    
    \item We extend the framework to noisy function evaluations with finite but \textit{a priori} unknown noise bounds. Noise bounds are estimated concurrently with the admissible gradient set, and we derive the theoretically optimal auxiliary sampling radius that minimises gradient set uncertainty. This results in a direct, noise-aware gradient estimation approach without smoothing, surrogate-based interpolation, or heuristic step-size rules.
\end{enumerate}

\noindent We illustrate the efficacy of our proposed method by comparing it with existing gradient estimators in two distinct tests: a large-scale numerical gradient accuracy benchmark, and an evaluations-based performance comparison in the context of a first-order optimisation method.

This paper is organised in seven sections. Section~\ref{sec:preliminaries} introduces the problem and the assumptions. Section~\ref{sec:results-gradient-estimation} discusses our main results on gradient estimation, along with discussions on its computational aspects on Section~\ref{sec:practical-aspects}. Our proposed Set-Based Adaptive Gradient Estimator is then introduced in Section~\ref{sec:sage}. Comparative statistical tests, results, and discussions are in Section~\ref{sec:benchmark-results}, and we synthesise our conclusions in Section~\ref{sec:conclusion}.

\subsubsection*{Notation}
We consider a scalar function $f : \mathcal{X} \rightarrow \mathbb{R}$, where $\mathcal{X}\subset \mathbb{R}^D$ is a compact set, and $D$ is the dimensionality. We refer to a point as $\bm{x} \in \mathcal{X}$, its gradient as $\bm{g}(\bm{x}) \in \mathbb{R}^D$, and the Hessian as $\bm{H}(\bm{x}) \in \mathbb{R}^{D \times D}$. After $n \in \mathbb{N}$ function evaluations of $f$, we are given a data set 

\[
    \itrpt{X}{n}{} \doteq \left\{ (\idxpt{x}{i}{}, \idxval{z}{i}{}), i = 1, \ldots, n\right\}
\]

\noindent composed of tuples $(\idxpt{x}{i}{}, \idxval{z}{i}{})$, with $\idxpt{x}{i}{} \in \mathcal{X}$ being the evaluated point, and $\idxval{z}{i}{} \in \mathbb{R}$ the (noiseless) sampled value.  We use the superscript indexing $\cdot^{(i)}$ to refer to an element of a countable set, and $\cdot^{\langle n \rangle}$ to describe a quantity or set that changes with $n$. For simplicity, we refer to a sample $(\idxpt{x}{i}{}, \idxval{z}{i}{})$ just by its evaluated point, e.g., $\idxpt{x}{i}{} \in \itrpt{X}{n}{}$. For any two samples $\idxpt{x}{i}{}, \idxpt{x}{j}{}$, $\bm{a}_{ij} \doteq \idxpt{x}{j}{}-\idxpt{x}{i}{}$ is the vector difference, and $\bm{u}_{ij} \doteq \frac{\bm{a}_{ij}}{\|\bm{a}_{ij}\|}$ is the unit vector pointing from $\idxpt{x}{i}{}$ to $\idxpt{x}{j}{}$. As shorthand, we denote $\idxpt{g}{i}{}$ as the true gradient at $\idxpt{x}{i}{}$, $\idxpt{H}{i}{}$ the corresponding true Hessian, and the scalar $\idxval{H}{i}{} > 0$ its (spectral) norm (note the boldface and regular type difference). We denote with $\|\cdot\|$ the 2-norm when operating on a vector, and the spectral norm when on a matrix.

\section{Preliminaries}
\label{sec:preliminaries}

We consider the problem of finding the gradient $\idxpt{g}{i}{}$ at a sampled point $\idxpt{x}{i}{} \in \itrpt{X}{n}{}$, from a data set generated by sampling $f$. In this paper, we do not know the explicit expression for $f$, nor its basis function composition. However, we do have assumptions on its regularity:

\begin{assumption}
\label{ass:l-hessian}
    The considered function $f \in \mathcal{C}^2$, and has a Lipschitz continuous Hessian

    \[
        \frac{\|\bm{H}(\bm{x}_1) - \bm{H}(\bm{x}_2)\|}{\|\bm{x}_1-\bm{x}_2\|} \leq \gamma_H
    \]

    \noindent $\forall \bm{x}_1, \bm{x}_2 \in \mathcal{X}$ (not necessarily samples), where $\mathcal{X}$ is a compact set. 
\end{assumption}

\noindent Another assumption is in order, pertaining to the acquisition of samples:

\begin{assumption}
\label{ass:noiseless-sampling}
    For any $\idxpt{x}{i}{} \in \mathcal{X}$, the corresponding function value is accessible by sampling:
    \[
        \idxval{z}{i}{} = f(\idxpt{x}{i}{}).
    \]
\end{assumption}

\noindent Whilst Assumption~\ref{ass:l-hessian} seems to be restrictive at first, we impose it only on the compact domain $\mathcal{X}$ over which the samples can be taken, and it is only needed for the theoretical development in Section~\ref{sec:results-gradient-estimation}. We nevertheless demonstrate in the benchmark results in Section~\ref{sec:benchmark-results} that our method still works competitively in practice, even in cases where this is not fulfilled. Assumption~\ref{ass:noiseless-sampling} simplifies the development of our basic results on gradient estimation; however, we also elegantly tackle the case of noisy evaluations (Section~\ref{subsec:noisy-evaluations}) using our proposed framework.

\section{Results on gradient estimation and refinement}
\label{sec:results-gradient-estimation}

We now give the first result regarding the accuracy of two-sample-based directional slope computation with respect to the actual directional derivative at $\idxpt{x}{i}{}$. This applies to any pair of samples, even far-away ones, arguing that, under Assumption \ref{ass:l-hessian}, global information from other samples can still provide local information on the gradient at the sample of interest.

\begin{corollary}[from Nesterov \cite{nesterov2018lectures}, Lemma 4.1.1]
\label{corollary:gradient-accuracy-2samples}
    Let Assumptions~\ref{ass:l-hessian}-\ref{ass:noiseless-sampling} hold. Given two samples $(\idxpt{x}{i}{},\idxval{z}{i}{}), (\idxpt{x}{j}{},\idxval{z}{j}{})$, it applies that $\idxpt{g}{i}{} \in \mathcal{G}_{ij}$ for which the following holds:

    \begin{equation}
        \mathcal{G}_{ij} \doteq \left\{\bm{g} :\Big| \frac{\idxval{z}{j}{} - \idxval{z}{i}{}}{\|\bm{a}_{ij}\|} - \bm{g}^\top \bm{u}_{ij} \Big| \leq \frac{1}{6} \gamma_H \|\bm{a}_{ij}\|^2 + \frac{1}{2} \idxval{H}{i}{} \|\bm{a}_{ij}\|\right\}.
    \label{eqn:gradient-slab}
    \end{equation}
\end{corollary}

{\it Proof:} See Appendix.

\begin{figure}
\centering
    \includegraphics[width=0.5\textwidth]{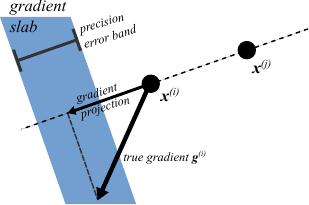}
    \caption{Visualisation of the gradient slab defined by \\points $\idxpt{x}{i}{}$ and $\idxpt{x}{j}{}$ in two dimensions.}
    \label{fig:gradient-slab}
\end{figure}

The result highlights that the worst-case error between the two-sample-based directional slope estimate $\frac{\idxval{z}{j}{} - \idxval{z}{i}{}}{\|\bm{a}_{ij}\|}$ and the projection of $\idxpt{g}{i}{}$ on $\bm{u}_{ij}$, worsens quadratically with the distance between the two points. Corollary~\ref{corollary:gradient-accuracy-2samples} also means that given only two samples $\idxpt{x}{i}{}$ and $\idxpt{x}{j}{}$, we can already extract information regarding $\idxpt{g}{i}{}$ through $\mathcal{G}_{ij}$, representing a set of admissible vectors (``gradient slab'') which we visualise in Fig.~\ref{fig:gradient-slab}. 


\begin{remark}[On gradient locality and sample informativeness]
Although the gradient is a local operator, Corollary~\ref{corollary:gradient-accuracy-2samples} highlights that information constraining the gradient at a query point need not originate exclusively from samples in its infinitesimal neighbourhood. The inequality in Corollary~\ref{corollary:gradient-accuracy-2samples} is induced by any pair of samples and does not rely on locality assumptions such as symmetric sampling, small perturbation radii, or geometric poisedness. As a result, samples that are not generated as dedicated auxiliary perturbations---potentially acquired at earlier iterations or for other purposes---can still provide valid constraints on the gradient.
  
This contrasts with most estimator-first approaches to gradient approximation, in which locality is enforced implicitly by the construction of the estimator (e.g., finite differences, simplex gradients, or smoothing-based estimators), and samples are generated solely for the purpose of estimating the gradient at the query point \cite{audet2017derivative,conn2009introduction,nesterov2017random}. From a set membership perspective, locality is not a prerequisite: what matters is whether a sample pair induces a valid inequality on the unknown gradient under the assumed smoothness conditions. This observation motivates the reuse of existing samples and the exploitation of global information for gradient set identification.
\end{remark}

\begin{remark}[Relation to existing smoothness-based bounds and novelty of the slab construction]
Corollary~\ref{corollary:gradient-accuracy-2samples} is a direct consequence of standard Taylor-remainder bounds for twice-differentiable functions with Lipschitz-continuous Hessian (see, e.g., \cite{nesterov2018lectures}). Related smoothness assumptions are ubiquitous in derivative-free optimisation, where they are typically used to derive \textit{a posteriori} error bounds for specific point estimators (finite differences, simplex gradients, interpolation models) under prescribed sampling geometries (e.g., \cite{audet2017derivative,conn2009introduction,hare2022error,regis2014calculus}).

Our methodological use of the bound is different: rather than bounding the error of a chosen estimator in a norm form, we rewrite the two-point remainder inequality as a \textit{distance-normalised directional-slope constraint} that is \textit{linear in the unknown gradient} at the query point. Each sample pair thus induces a pair of halfspace constraints (a ``slab'') on $\idxpt{g}{i}{}$, and the intersection of these slabs yields a polyhedral admissible gradient set whose geometry (e.g., diameter, unboundedness) directly quantifies uncertainty and drives adaptive sampling.

This slab construction relies on having an explicit uniform bound on the second-order Taylor remainder along the segment connecting two samples. The Lipschitz-Hessian assumption provides such a bound with controlled growth in the inter-sample distance (quadratic in $\|\bm{a}_{ij}\|$), which makes nonlocal sample pairs informative within the set membership framework. Under weaker smoothness assumptions, analogous constraints may still be derived, but they typically inflate faster with distance, reducing the practical utility of nonlocal slabs and weakening refinement guarantees.
\end{remark}

For brevity in the subsequent discussion, we now denote the two-sample-based directional slope estimate $\frac{\idxval{z}{j}{} - \idxval{z}{i}{}}{\|\bm{a}_{ij}\|}$ as $\tilde{s}_{ij}$. As a consequence of Corollary~\ref{corollary:gradient-accuracy-2samples}, we state the following result:

\begin{theorem}
\label{theorem:gradient-set}
    Let Assumptions~\ref{ass:l-hessian}-\ref{ass:noiseless-sampling} hold, and consider a point $\idxpt{x}{i}{} \in \itrpt{X}{n}{}$. Given $\idxval{H}{i}{}$ and $\gamma_H$, and the data set $\itrpt{X}{n}{}$, the gradient $\idxpt{g}{i}{}$ at $\idxpt{x}{i}{}$ is guaranteed to belong to the set

    \begin{equation}
        \mathcal{G} \doteq 
        \left\{
        \bm{g} ~:~
        \begin{bmatrix}
            -\bm{u}_{ij}^\top \\
             \bm{u}_{ij}^\top
        \end{bmatrix}
        \bm{g}
        \leq
        \begin{bmatrix}
            -\tilde{s}_{ij} \\ 
             \tilde{s}_{ij}
        \end{bmatrix}
        +
        \frac{1}{2}
        \begin{bmatrix}
            \|\bm{a}_{ij}\| \\ 
            \|\bm{a}_{ij}\|
        \end{bmatrix}
        \idxval{H}{i}{}
        +
        \frac{1}{6}
        \begin{bmatrix}
            \|\bm{a}_{ij}\|^2 \\ 
            \|\bm{a}_{ij}\|^2 \\
        \end{bmatrix}
        \gamma_H, \forall \idxpt{x}{j}{} \in \itrpt{X}{n}{} \setminus \idxpt{x}{i}{}
        \right\}
    \label{eqn:gradient-set}
    \end{equation}    
\end{theorem}

{\it Proof:} See Appendix.

\begin{figure}[b]
\centering
    \includegraphics[width=0.5\textwidth]{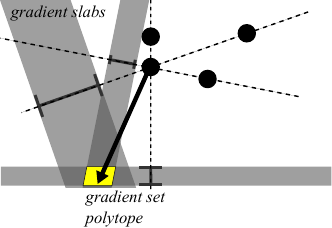}
    \caption{Visualisation of the gradient set polytope in two dimensions.}
    \label{fig:gradient-polytope}
\end{figure}

The set \eqref{eqn:gradient-set} is a convex polytope given by the intersection of $n-1$ slabs, as shown in Fig.~\ref{fig:gradient-polytope}.

\subsection{\it Estimating the gradient set}

Our results in Theorem~\ref{theorem:gradient-set} are given however with restrictive (and in most cases, unrealistic) assumptions of known $\idxval{H}{i}{}$ (Hessian spectral norm) and $\gamma_H$. Therefore, in building our proposed method, we resort to finding their corresponding estimates $\idxval{\tilde{H}}{i}{}$ and $\tilde{\gamma}_H$, whilst finding the tightest unfalsified set of gradients $\tidxpt{g}{i}{}$ in a set membership sense. Moving the estimated quantities in \eqref{eqn:gradient-set} to the left side, we have

\begin{equation}
    \underbrace{\begin{bmatrix}
        -\bm{u}_{ij_1}^\top & -\frac{1}{2}\|\bm{a}_{ij_1}\| & -\frac{1}{6}\|\bm{a}_{ij_1}\|^2 \\
         \bm{u}_{ij_1}^\top & -\frac{1}{2}\|\bm{a}_{ij_1}\| & -\frac{1}{6}\|\bm{a}_{ij_1}\|^2 \\
        -\bm{u}_{ij_2}^\top & -\frac{1}{2}\|\bm{a}_{ij_2}\| & -\frac{1}{6}\|\bm{a}_{ij_2}\|^2 \\
         \bm{u}_{ij_2}^\top & -\frac{1}{2}\|\bm{a}_{ij_2}\| & -\frac{1}{6}\|\bm{a}_{ij_2}\|^2 \\
         \vdots & \ddots & \vdots \\
        -\bm{u}_{ij_{n-1}}^\top & -\frac{1}{2}\|\bm{a}_{ij_{n-1}}\| & -\frac{1}{6}\|\bm{a}_{ij_{n-1}}\|^2 \\
         \bm{u}_{ij_{n-1}}^\top & -\frac{1}{2}\|\bm{a}_{ij_{n-1}}\| & -\frac{1}{6}\|\bm{a}_{ij_{n-1}}\|^2 \\
    \end{bmatrix}}_{\bm{A}}
    \underbrace{\begin{bmatrix}
        \tilde{\bm{g}}^{(i)} \\ \idxval{\tilde{H}}{i}{} \\ \tilde{\gamma}_H
    \end{bmatrix}}_{\tilde{\bm{v}}}
    \leq
    \underbrace{\begin{bmatrix}
        -\tilde{s}_{ij_1} \\ \tilde{s}_{ij_1} \\ -\tilde{s}_{ij_2} \\ \tilde{s}_{ij_2} \\ \vdots \\ -\tilde{s}_{ij_{n-1}} \\ \tilde{s}_{ij_{n-1}}
    \end{bmatrix}}_{\bm{b}}.
\end{equation}

We now propose to solve the following linear programme (LP) to find a feasible set of possible gradients using the data set:

\begin{align}
    \label{eqn:gradient-est-sage}
    \min_{\tilde{\bm{v}}} &
    \begin{bmatrix}
        0 & \ldots & 0 & 1 & 1
    \end{bmatrix}
    \tilde{\bm{v}} \\
    \label{eqn:gradient-est-sage2}
    \text{s.t.} & 
    \left[\begin{array}{c}
        \bm{A} \\
        \midrule
        \begin{matrix}
        0 & \cdots & 0 & -1 & 0 \\
        0 & \cdots & 0 & 0 & -1
        \end{matrix}
    \end{array}\right] \tilde{\bm{v}} \leq 
    \left[\begin{array}{c}
        \bm{b} \\
        \midrule
        0 \\
        0        
    \end{array}\right]
\end{align}

\noindent The above programme minimises the sum of the (non-negative) estimates $\idxval{\tilde{H}}{i}{}, \tilde{\gamma}_H$ such that the polytope is non-empty, which translates to finding the tightest estimated gradient set $\tilde{\mathcal{G}} \doteq \bigcap_j \tilde{\mathcal{G}}_{ij}$. Note that $\bm{A}$ is composed of vector differences of the other samples with respect to the point of interest and their norms, whilst $\bm{b}$ is built from the respective directional slope quantities. At the same time, we also achieve a gradient estimate $\tilde{\bm{g}}^{(i)}$, such that:

\begin{equation}
\label{eqn:gradient-set-estimate}
    \tilde{\bm{g}}^{(i)} \in \tilde{\mathcal{G}} \doteq \Big\{\bm{g} \in \mathbb{R}^D : \bm{A}_l \bm{g} \leq \underbrace{\bm{b} - \bm{A}_r \left[\idxval{\tilde{H}}{i}{} ~ \tilde{\gamma}_H\right]^\top}_{\bm{b}'} \Big\}
\end{equation}

\noindent where $\bm{A}_l$ is the left $D$ columns of $\bm{A}$, and $\bm{A}_r$ is composed of its other columns (note the correspondence with \eqref{eqn:gradient-set}).

\subsection{\it Uncertainty quantification and point gradient estimate selection}
To quantify the precision of the gradient set $\tilde{\mathcal{G}}$, it is useful to quantify the worst-case uncertainty in choosing an estimate $\tidxpt{g}{i}{}$, as well as the directions along with the uncertainty is highest. We define the \textit{gradient set diameter} as follows:

\begin{definition}
\label{defn:gradient-set-diameter}
    The diameter of a gradient set $\tilde{\mathcal{G}}$ is given by

    \begin{equation}
        \rho(\tilde{\mathcal{G}}) \doteq \max_{\bm{g}_1, \bm{g}_2 \in \tilde{\mathcal{G}}} \left\|\bm{g}_1 - \bm{g}_2\right\|.
    \end{equation}
    
\end{definition}

Correspondingly, the direction of maximum uncertainty is $\hat{\bm{e}} = \frac{\bm{g}_1 - \bm{g}_2}{\|\bm{g}_1 - \bm{g}_2\|}$. Recognising that computing $\rho(\tilde{\mathcal{G}})$ is a nonconvex problem, we efficiently address this using the centroid of a rotated bounding hyperbox, described as follows.

Recalling $\bm{u}_{ij} = \frac{\bm{a}_{ij}}{\|\bm{a}_{ij}\|}$, for every selected neighbour $\idxpt{x}{j}{}$, we stack all these unit vectors as 

\[
    \bm{U} = \left[
    \begin{matrix}
        \bm{u}_{ij_1}^\top \\
        \bm{u}_{ij_2}^\top \\
        \vdots \\
        \bm{u}_{ij_{n-1}}^\top \\
    \end{matrix}
    \right].
\]

\noindent We then compute the singular value decomposition (SVD) of $\bm{U}$ to obtain $\bm{Q}$, with rows $\bm{q}_k$ being an orthonormal basis for $\mathbb{R}^D$ with directions most to least covered by $\bm{u}_{ij}$. As these $\bm{u}_{ij}$ are also the constraint normals of $\tilde{\mathcal{G}}$ (see the slab construction in Fig.~\ref{fig:gradient-slab}), $\bm{Q}$ contains a natural basis for our rotated hyperbox. On the other hand, when $\bm{U}$ is rank-deficient, the corresponding $\bm{q}_k$ is in the null space of $\bm{U}$. 

For each $\bm{q}_k$ we solve two linear programmes

\begin{align}
    \underline{w}_k &= \min_{\bm{g} \in \tilde{\mathcal{G}}} \bm{q}_k^\top \bm{g} \label{eqn:w_min}\\
    \overline{w}_k &= \max_{\bm{g} \in \tilde{\mathcal{G}}} \bm{q}_k^\top \bm{g} \label{eqn:w_max}
\end{align}

\noindent which may return $\underline{w}_k = -\infty$ and/or $\overline{w}_k = +\infty$ if $\tilde{\mathcal{G}}$ is unbounded along $\bm{q}_k$. With $w_k = \overline{w}_k - \underline{w}_k$ and stacking them to $\bm{w}$, the \textit{approximate} diameter is

\begin{equation}
\label{eqn:gradient-set-diameter}
    \tilde{\rho}(\tilde{\mathcal{G}}) = \| \bm{w}\|,
\end{equation}

\noindent the diagonal length of the bounding hyperbox, which upper-bounds the true diameter since $\tilde{\mathcal{G}}$ is contained in the box. The direction of maximum uncertainty $\hat{\bm{e}}$ is then given by:

\begin{equation}
\label{eqn:direction-maximum-uncertainty}
    \hat{\bm{e}} =  
    \begin{cases}
        \frac{\bm{Q}^\top \bm{w}}{\|\bm{w}\|} & \text{if $\tilde{\rho}(\tilde{\mathcal{G}}) < +\infty$ (the unit-norm diagonal direction)} \\
        \bm{q}_{k^*} & \text{if $\tilde{\rho}(\tilde{\mathcal{G}}) = +\infty$ for any $k^*$ with $\overline{w}_{k^*} = +\infty$ (or $\underline{w}_{k^*}=-\infty$)}
    \end{cases}
\end{equation}

\noindent If the second case is due to $\bm{U}$ being rank-deficient, $\hat{\bm{e}}$ is perpendicular to all $\bm{u}_{ij}$ by construction of the SVD.

Furthermore, in most applications, a point gradient estimate $\tidxpt{g}{i}{}$ is required rather than $\tilde{\mathcal{G}}$ itself. Conveniently, we can use the same hyperbox-based infrastructure for solving $\tilde{\rho}$ to pick a sensible $\tidxpt{g}{i}{}$. If $\tilde{\rho}(\tilde{\mathcal{G}}) < +\infty$, we can consider $m_k = \frac{1}{2}(\overline{w}_k + \underline{w}_k)$ and stack them to build $\bm{m}$; the point estimate is now 

\begin{equation}
\label{eqn:gradient-estimate}
    \tidxpt{g}{i}{} = \bm{Q}^\top \bm{m}.
\end{equation}

\subsection{\it Refining the estimated gradient set by strategic sampling}
\label{subsec:strategic-sampling}

Depending on the application, there can be requirements on the gradient estimate precision, quantified in this paper as a desired gradient set diameter $\rho^*$. Building on our results, we propose a simple iterative procedure to refine the gradient estimate.

\begin{enumerate}
    \item Compute the (approximate) diameter $\tilde{\rho}(\tilde{\mathcal{G}})$, using \eqref{eqn:w_min}-\eqref{eqn:gradient-set-diameter}.
    \item If $\tilde{\rho}(\tilde{\mathcal{G}}) < \rho^*$, then we already have the desired gradient accuracy, and we use the gradient estimate $\tilde{\bm{g}}^{(i)}$. Otherwise, compute the direction of maximum uncertainty $\hat{\bm{e}}$ using \eqref{eqn:direction-maximum-uncertainty}.
    
    \item Generate an auxiliary sample at
    
    \begin{equation}
    \label{eqn:auxiliary-sample}
        \idxpt{x}{j}{} = \idxpt{x}{i}{} + \alpha \hat{\bm{e}}, 
    \end{equation}
    
    \noindent where $\alpha > 0$ is a small (or otherwise prescribed) positive value.
    
    \item Recompute $\tilde{\mathcal{G}}$ using \eqref{eqn:gradient-est-sage}-\eqref{eqn:gradient-est-sage2} and \eqref{eqn:gradient-set-estimate}, with $\itrpt{X}{n}{} \cup (\idxpt{x}{j}{}, \idxval{z}{j}{})$
    \item Go to Step 1.
\end{enumerate}

\subsection{\it Results with finite additive noise}
\label{subsec:noisy-evaluations}

Conventional gradient estimation techniques, e.g., finite-difference-based methods, rely on accurate function evaluations, and are therefore highly sensitive to noise. In industrial applications, e.g., in experiment-based optimisation, we are faced with noisy samples (from measurement noise, process disturbance, or even numerical precision) that must be handled in gradient estimation. Here, we consider a setting where function evaluations are subject to noise with finite bounds,

\begin{equation}
\label{eqn:evaluation-with-noise}
    \idxpt{z}{i}{} = f(\idxpt{x}{i}{}) + \epsilon_i, |\epsilon_i| \leq \overline{\epsilon}.
\end{equation}

Our result is a slight modification to Corollary~\ref{corollary:gradient-accuracy-2samples}, which we now state:

\begin{corollary}[from Nesterov \cite{nesterov2018lectures}, Lemma 4.1.1]
\label{corollary:gradient-accuracy-2samples-noise}
    Given two samples $(\idxpt{x}{i}{}, \idxval{z}{i}{})$ and $(\idxpt{x}{j}{}, \idxval{z}{j}{})$, for which $\idxval{z}{i}{}$ and $\idxval{z}{j}{}$ are acquired in the manner of \eqref{eqn:evaluation-with-noise}, we have
    
    \begin{equation}
        \Big| \frac{\idxval{z}{j}{} - \idxval{z}{i}{}}{\|\bm{a}_{ij}\|} - \idxptT{g}{i}{} \bm{u}_{ij} \Big| \leq \frac{1}{2} \|\bm{a}_{ij}\| \idxval{H}{i}{} + \frac{1}{6}\|\bm{a}_{ij}\|^2 \gamma_H + \frac{2\overline{\epsilon}}{\|\bm{a}_{ij}\|}.
        \label{eqn:directional-derivative-uncertainty-noise}
    \end{equation}
\end{corollary}

{\it Proof:} See Appendix.

\noindent We note two important observations:

\begin{enumerate}
    \item There are now two sources for increased uncertainty with respect to the directional derivative estimate. When $\idxpt{x}{i}{}$ and $\idxpt{x}{j}{}$ are too far apart, the first two terms in the right-hand side of \eqref{eqn:directional-derivative-uncertainty-noise} dominate (corresponding to that of the noiseless result). Now in addition, when the samples are too close, the last (additional) term dominates, resulting in the estimates ``drowning in noise''.
    \item We cannot simply sample as close as possible to $\idxpt{x}{i}{}$ to refine the gradient set to a desired diameter. In fact, due to noise, there are cases where our desired diameter is impossible to achieve due to high noise bounds.
\end{enumerate}

\begin{remark}
Corollary~\ref{corollary:gradient-accuracy-2samples-noise} extends the inequality of Corollary~\ref{corollary:gradient-accuracy-2samples} to the case of bounded additive noise, whilst preserving the same set membership interpretation and geometric role. Error bounds for noisy finite-difference or interpolation-based gradient estimators have been previously studied under bounded noise \cite{berahas2019linear,more2012estimating,scheinberg2022finite}, to analyse the accuracy or optimal step size of a point estimator. In contrast, Corollary~\ref{corollary:gradient-accuracy-2samples-noise} keeps the two-sample directional formulation of Corollary~\ref{corollary:gradient-accuracy-2samples} and incorporates noise directly into the slab constraints, enabling the joint treatment of gradient uncertainty, curvature effects, and noise within the same admissible-set construction. As in Corollary~\ref{corollary:gradient-accuracy-2samples}, the role of the inequality is not to certify a particular estimator, but to define a feasible region guaranteed to contain the true gradient. The geometry of this feasible region (e.g., boundedness or diameter) is then exploited for data-driven refinement and noise-aware gradient selection.
\end{remark}

Taking the minimum value of the right-hand side of \eqref{eqn:directional-derivative-uncertainty-noise}, the optimal auxiliary sample distance $\alpha^*$ used in \eqref{eqn:auxiliary-sample} is given by

\begin{equation}
    \alpha^* =  \min \left\{ \alpha \in \mathbb{R}_{>0} : \frac{1}{3} \gamma_H \alpha^3 + \frac{1}{2} \idxval{H}{i}{} \alpha^2 - 2 \overline{\epsilon} = 0 \right\},
    \label{eqn:optimal-auxiliary-spacing-noise}
\end{equation}

\noindent and substituting back $\alpha^*$ to the right hand side of \eqref{eqn:directional-derivative-uncertainty-noise}, we obtain the \textit{theoretical best achievable} gradient estimate precision $\rho^*(\tilde{\mathcal{G}})$. When $\overline{\epsilon} = 0$, \eqref{eqn:optimal-auxiliary-spacing-noise} collapses into $\alpha^* = 0$ and $\rho^*(\tilde{\mathcal{G}}) = 0$, which is why we can set $\alpha$ to any small finite number in \eqref{eqn:auxiliary-sample} in the noiseless case (in this paper, we set $\alpha = 1\times 10^{-6}$ for the noiseless case as the nominal auxiliary step size).

\begin{remark}
    Equation \eqref{eqn:optimal-auxiliary-spacing-noise} is a theoretical optimal sampling radius given the \textit{a priori} knowledge of $\idxval{H}{i}{}$ and $\gamma_H$. In a real use case, we can fall back to using their respective estimates, which, whilst losing the theoretical guarantees, do still work in practice, as shown later in the benchmark tests.
\end{remark}

Considering our desired gradient set diameter $\rho$, we modify the condition in step~2 in Section~\ref{subsec:strategic-sampling} to be $\tilde{\rho}(\tilde{\mathcal{G}}) < \max(\rho^*, \rho^*(\tilde{\mathcal{G}}))$. This means that with noisy evaluations, we should refine our gradient set diameter until either our desired precision, or the theoretical best one is achieved.

In the case that $\overline{\epsilon}$ is not \textit{a priori} known, we can also use the data set $\itrpt{X}{n}{}$ to directly generate an estimate $\tilde{\epsilon}$, together with $\tilde{\bm{g}}^{(i)}, \tilde{H}^{(i)}, \tilde{\gamma}_H$. Using \eqref{eqn:directional-derivative-uncertainty-noise}, equations \eqref{eqn:gradient-est-sage}-\eqref{eqn:gradient-est-sage2} become

\begin{align}
    \label{eqn:gradient-est-sage-noise}
    \min_{\bm{v}} &
    \begin{bmatrix}
        0 & \ldots & 0 & 1 & 1 & 1
    \end{bmatrix}
    \tilde{\bm{v}}_\epsilon \\
    \label{eqn:gradient-est-sage-noise2}
    \text{s.t.} & 
    \left[\begin{array}{c}
        \bm{A}_\epsilon \\
        \midrule
        \begin{matrix}
        0 & \cdots & 0 & -1 & 0 & 0 \\
        0 & \cdots & 0 & 0 & -1 & 0 \\
        0 & \cdots & 0 & 0 & 0 & -1        
        \end{matrix}
    \end{array}\right] \tilde{\bm{v}}_\epsilon \leq 
    \left[\begin{array}{c}
        \bm{b} \\
        \midrule
        0 \\
        0 \\
        0
    \end{array}\right]
\end{align}

\noindent with

\[
    \bm{A}_\epsilon = 
    \begin{bmatrix}
        -\bm{u}_{ij_1}^\top & -\frac{1}{2}\|\bm{a}_{ij_1}\| & -\frac{1}{6}\|\bm{a}_{ij_1}\|^2  & -\frac{2}{\|\bm{a}_{ij_1}\|} \\
         \bm{u}_{ij_1}^\top & -\frac{1}{2}\|\bm{a}_{ij_1}\| & -\frac{1}{6}\|\bm{a}_{ij_1}\|^2  & -\frac{2}{\|\bm{a}_{ij_1}\|} \\
         \vdots & \ddots & \vdots & \vdots \\
        -\bm{u}_{ij_{n-1}}^\top & -\frac{1}{2}\|\bm{a}_{ij_{n-1}}\| & -\frac{1}{6}\|\bm{a}_{ij_{n-1}}\|^2  & -\frac{2}{\|\bm{a}_{ij_{n-1}}\|} \\
         \bm{u}_{ij_{n-1}}^\top & -\frac{1}{2}\|\bm{a}_{ij_{n-1}}\| & -\frac{1}{6}\|\bm{a}_{ij_{n-1}}\|^2  & -\frac{2}{\|\bm{a}_{ij_{n-1}}\|} \\
    \end{bmatrix}
\]

\noindent and 

\begin{equation}
    \tilde{\bm{v}}_\epsilon = [~\tilde{\bm{g}}^{(i)\top} ~ \idxval{\tilde{H}}{i}{} ~ \tilde{\gamma}_H ~ \tilde{\epsilon}~]^\top.
    \label{eqn:noise-estimation-variable-definition}
\end{equation}

\noindent The gradient set is correspondingly

\begin{equation}
\label{eqn:gradient-set-estimate-noisy}
    \tilde{\mathcal{G}} = \Big\{\bm{g} \in \mathbb{R}^D : \bm{A}_{\epsilon l} \bm{g} \leq \underbrace{\bm{b} - \bm{A}_{\epsilon r} \left[\idxval{\tilde{H}}{i}{} ~ \tilde{\gamma}_H ~ \tilde{\epsilon} \right]^\top}_{\bm{b}_\epsilon'} \Big\}
\end{equation}

\noindent where $\bm{A}_{\epsilon l}$ is the left column of $\bm{A}_\epsilon$, and $\bm{A}_{\epsilon r}$ is its three rightmost columns. With these results we now have a \textit{systematic} approach on deciding how far from the query point we should sample, to have the greatest information gain regarding $\idxpt{g}{i}{}$. Furthermore, even with \textit{a priori} unknown noise bound $\overline{\epsilon}$, we have an LP-based approach to estimate it simultaneously with the other quantities, seamlessly integrated in our proposed set membership-based framework.

\begin{remark}[On auxiliary sampling radius and step-size heuristics]
Corollary~\ref{corollary:gradient-accuracy-2samples-noise} makes explicit the competing effects of curvature magnitude, curvature variation, and evaluation noise on the uncertainty induced by a pair of samples. In estimator-first approaches, this trade-off is typically resolved by selecting an auxiliary sampling radius (or step size) that minimises the error of a point-valued directional derivative or gradient estimate. In practice, such step sizes are often chosen using heuristics, problem-dependent tuning, or rules of thumb, and may require retuning when noise characteristics or problem scales change \cite{boresta2022mixed,more2012estimating,scheinberg2022finite}.
  
In contrast, within a set membership framework the sampling radius directly affects the geometry of the admissible gradient set. Corollary~\ref{corollary:gradient-accuracy-2samples-noise} provides a principled way to select the sampling distance that minimises the worst-case uncertainty of the gradient set (ultimately arriving at \eqref{eqn:optimal-auxiliary-spacing-noise}), rather than the error of a specific estimator. This shifts the role of auxiliary sampling from estimator calibration to uncertainty reduction: additional samples are acquired only when they are expected to meaningfully tighten the admissible set. As a consequence, step-size selection emerges from the uncertainty structure itself, rather than being imposed externally.
\end{remark}

\section{Practical implementation and computational aspects}
\label{sec:practical-aspects}

\subsection{\it Bootstrapping and adaptive sampling scale}

At a new query point $\idxpt{x}{i}{}$, and especially without an initial samples set, the admissible gradient set $\tilde{\mathcal{G}}$ is unbounded due to the absence of sufficient constraints. In this case, we implement bootstrapping: SAGE first constructs a symmetric coordinate-aligned stencil:

\begin{equation}
    \idxpt{x}{i}{} \pm \alpha_o \hat{\bm{e}}_d, \quad d = 1, \dots, D,
\end{equation}

\noindent where $\alpha_o > 0$ is a small initial step size and $\hat{\bm{e}}_d$ is a unit vector at coordinate direction $d$. This produces $2D$ samples, generating paired gradient slabs that bound each component of the gradient. In the noiseless case ($\overline{\epsilon} = 0.0$), the initial samples already yield highly accurate estimates, and there is no need for further auxiliary sampling-based gradient set refinement. Otherwise, in the presence of noise, small step sizes lead to noise-dominated constraints (see Corollary~\ref{corollary:gradient-accuracy-2samples-noise}) resulting in an inflated admissible set. Nevertheless, this initial stencil remains informative, as the resulting LP provides a noise bound estimate $\tilde{\epsilon}$.

Whilst \eqref{eqn:optimal-auxiliary-spacing-noise} computes the best auxiliary sampling radius, it may not give a positive root, which can occur in early stages when $\tilde{H}^{(i)} \approx 0$ and $\tilde{\gamma}_H \approx 0$. In such cases, a bootstrap step size $\alpha = 2\sqrt{\tilde{\epsilon}}$ is used, providing a conservative estimate of the noise-appropriate sampling scale. Accordingly, a subsequent symmetric stencil is constructed: $\idxpt{x}{i}{} \pm \alpha \hat{\bm{e}}_d$, yielding well-conditioned constraints that significantly reduce the diameter of $\tilde{\mathcal{G}}$ and ensure all coordinate directions are properly bounded. Subsequent refinement steps operate using samples at this informative scale.

\subsection{\it Constraints filtering for computational efficiency}

The computational complexity of the LP \eqref{eqn:gradient-est-sage}-\eqref{eqn:gradient-est-sage2} depends on the underlying algorithm used, e.g., simplex method, interior point method, or ellipsoid method. Therefore, computing $\tilde{\mathcal{G}}$ can have complexity ranging from polynomial to exponential with respect to constraints (i.e., the rows of $\bm{A}$). Furthermore, a neighbouring samples set, whose displacement vectors from $\idxpt{x}{i}{}$ are dominated by near-parallel directions, can cause numerical degeneracy in the LP as well. To reduce the computational complexity whilst preserving the richness of directions, we select a subset of neighbours using a combination of column-pivoted QR and proximity ranking to $\idxpt{x}{i}{}$. 

We therefore select $N_f$ neighbours (which we set in this paper as $N_f=\min(5D, n)$) using a column-pivoted QR on their respective unit directions $\bm{u}_{ij}$, presorted by $\left| \|\idxpt{x}{j}{} - \idxpt{x}{i}{}\| - \alpha^* \right|$ to break ties. This mechanism greedily picks, at each step, the remaining direction most linearly independent from those already chosen (the standard rank-revealing criterion) \cite{golub2013matrix}. The radius therefore only breaks ties among already-redundant directions, whilst linear independence primarily drives the neighbour selection.

\subsection{\it Stabilisation and termination}

The estimation of latent quantities $\idxval{\tilde{H}}{i}{}$ and $\tilde{\gamma}_H$ may fluctuate as new samples are added. To mitigate this, monotonically increasing updates are enforced. As  $\idxval{\tilde{H}}{i}{}$ is associated to a single query point, monotone updates are held only for the same $\idxpt{x}{i}{}$, and is reset to a small default value when we estimate the gradient at another query point (e.g., across iterates for a first-order method). On the other hand, $\tilde{\gamma}_H$ holds across the entire $\mathcal{X}$, and this estimate is therefore preserved and updated monotonically across all query points. 

Furthermore, gradient refinement is terminated when successive updates yield negligible changes in the gradient estimate, indicating convergence of the admissible set. A maximum number of auxiliary evaluations is also imposed as a safeguard against pathological cases.

\section{Set-Based Adaptive Gradient Estimator (SAGE)}
\label{sec:sage}

Based on our results, we are now in a position to introduce our proposed method Set-Based Adaptive Gradient Estimator (SAGE). We give here a high-level summary of our algorithm in pseudo-code form (see Algorithm \ref{algo:sage}), for which a Python-based implementation is publicly available at \url{https://github.com/lorenzosabugjr/SAGE}.

\begin{algorithm}
    \DontPrintSemicolon
    \KwData{samples set $\itrpt{X}{n}{}$, query point $\idxpt{x}{i}{}$}    
    Compute $\tilde{\mathcal{G}}$ using \eqref{eqn:gradient-est-sage}-\eqref{eqn:gradient-est-sage2} and \eqref{eqn:gradient-set-estimate} (or \eqref{eqn:gradient-est-sage-noise}-\eqref{eqn:gradient-est-sage-noise2} and \eqref{eqn:gradient-set-estimate-noisy})\label{algo-start}\;
    Estimate gradient set diameter $\tilde{\rho}(\tilde{\mathcal{G}})$ with \eqref{eqn:gradient-set-diameter}\;
    \If{gradient set diameter $\tilde{\rho}(\tilde{\mathcal{G}}) > \max(\rho^*, \rho^*(\tilde{\mathcal{G}}))$}
    {
        Compute direction of maximum uncertainty $\hat{\bm{e}}$ using \eqref{eqn:direction-maximum-uncertainty}\;
        Get next sampling point $\idxpt{x}{n+1}{}$ using \eqref{eqn:auxiliary-sample} \;
        Evaluate $\idxpt{x}{n+1}{}$, update data set $\itrpt{X}{n+1}{} \leftarrow \itrpt{X}{n}{} \cup (\idxpt{x}{n+1}{}, \idxval{z}{n+1}{})$\;
        Update number of evaluations $n \leftarrow n+1$\;
        Go to Line~\ref{algo-start}\;
    }
    Compute point gradient estimate using \eqref{eqn:w_min}-\eqref{eqn:gradient-estimate} \;
    \KwResult{estimated gradient $\tilde{\bm{g}}^{(i)}$}
\caption{Set-Based Adaptive Gradient Estimator (SAGE)}
\label{algo:sage}
\end{algorithm}

\section{Benchmark and comparative results}
\label{sec:benchmark-results}

We compared SAGE with five different gradient estimators from the literature and common practice, to assess the reliability of gradient estimators in different use contexts: as standalone gradient estimators, and as a part of a gradient descent-based optimisation method. A summary of the compared gradient estimators, along with their hyperparameters as found in the literature, is given in Table~\ref{table:benchmark-gradient-estimators}. The first two were selected as the most commonly used gradient estimators in practice. GSG (Gaussian Smoothing Gradient) and CGSG (Central Gaussian Smoothing Gradient) enjoy increased popularity in the machine learning community because of their simplicity and their theoretical convergence properties. Lastly, NMXFD (Normalized Mixed Finite Difference) is a recently-proposed method for which an open-source code is available.

\begin{table}[]
\small
\centering
\begin{tabular}{lllp{2.0in}}
\toprule
\multirow{2}{*}{\textbf{Name}} & \multicolumn{2}{c}{\textbf{Auxiliary samples}}                        & \multirow{2}{*}{\textbf{Notes}}                                                                                     \\ \cline{2-3}
                               & \multicolumn{1}{l}{Number} & Nominal step size       &                                                                                                                     \\ \midrule
FFD \cite{conn2009introduction}                            & \multicolumn{1}{l}{$D$}    & $1\times10^{-6}$ &                                                                                                                    \\ 
CFD \cite{conn2009introduction}                           & \multicolumn{1}{l}{$2D$}   & $1\times10^{-6}$ &                                                                                                                    \\
GSG \cite{berahas2022comparison}                           & \multicolumn{1}{l}{$D$}    & $1\times10^{-6}$ & Given as ``smoothing radius''                                                                                                                    \\ 
CGSG \cite{berahas2022comparison}                          & \multicolumn{1}{l}{$2D$}   & $1\times10^{-6}$ & Given as ``smoothing radius''                                                                                                                    \\
NMXFD \cite{boresta2022mixed}                         & \multicolumn{1}{l}{$4D$}   & $1\times10^{-2}$ & Given as ``scale parameter'' $\varsigma$. For each coordinate, 4 equally-spaced auxiliary samples in interval [-2.0, +2.0]$\varsigma$ from $\idxpt{x}{i}{}$. \\
SAGE (ours)                           & Adaptive & Adaptive & Initial auxiliary samples with step size $1\times10^{-6}$ ($2D$ in forward-/backward directions, like CFD), afterwards adaptive.           \\ \bottomrule
\end{tabular}
\caption{Summary of compared gradient estimators}
\label{table:benchmark-gradient-estimators}
\end{table}

\begin{table}[]
\centering
\small
\begin{tabular}{lp{1.0in}lp{0.8in}}
    \toprule
~ & \textbf{Name} & \textbf{Definition} $f(\bm{x})$ & \textbf{Notes} \\
    \midrule
    \textbf{P1} & Simple least squares (LSQ) & $\frac{1}{2} \|\bm{y} - \bm{Q}\bm{x}\|^2$                    &                \\
    \textbf{P2} & $L_1$-regularised least squares (``LASSO'') &  $\frac{1}{2} \|\bm{y} - \bm{Q}\bm{x}\|^2 + \lambda\|\bm{x}\|_1$ & Assumption~\ref{ass:l-hessian} does not hold \\ 
    \textbf{P3} & log-sum-exp &  $\log\left( \sum_{i=1}^D \exp\left( (\bm{Q}\bm{x})_i - y_i \right) \right) + \frac{\lambda}{2} \|\bm{x}\|^2$                   &                \\ 
    \textbf{P4} & $L_1$-regularised logistic regression (L1-log-reg) & $\log\left(1+\exp(-\bm{y}^\top\bm{Q}\bm{x})\right) + \lambda\|\bm{x}\|_1$ & Assumption~\ref{ass:l-hessian} does not hold \\ 
    \textbf{P5} & $L_2$-regularised logistic regression (L2-log-reg) & $\log\left(1+\exp(-\bm{y}^\top\bm{Q}\bm{x})\right) + \frac{\lambda}{2}\|\bm{x}\|^2$                    &                \\ 
    \bottomrule
\end{tabular}
\caption{Benchmark functions considered in this paper ($D=20$)}
\label{table:benchmark-problems}
\end{table}

\subsection{\it Benchmark functions}

In our tests, we chose five different benchmark functions (all with $D=20$), given in Table~\ref{table:benchmark-problems}. Whilst P1, P3, and P5 satisfy Assumption~\ref{ass:l-hessian}, functions P2 and P4 do not. Therefore, P2 and P4 are included as a stress test to assess practical robustness beyond our theoretical assumptions.

For each function P1-P5, we have considered 100 random function instances, by generating random positive definite matrices $\bm{Q}$ to eliminate as much as possible function structure-based bias that could favour any one gradient estimator.

\subsection{\it Gradient accuracy benchmark tests}
We tested the efficacy of each compared estimator by using it to query a gradient estimate at discrete points. For each of the 100 instances of P1-P5, we randomly generated 250 query points, making a grand total of 25000 gradient queries for each benchmark function. 

To ensure fairness across gradient estimators, all random generations used the same set of seeds. Furthermore, we reset/cleared the SAGE samples set for each new gradient query, and did \textit{not} reuse samples from previous queries on the same problem instance. All tests were repeated with different (uniform) noise bounds $\overline{\epsilon} = \{0.0, 1\times10^{-3}, 1.0\}$.

\subsubsection*{\it Comparison metric}
This benchmark considers the relative gradient accuracy $\delta_r > 0$ as the main metric 
    
\[
    \delta_r = \frac{\|\tilde{\bm{g}}(\idxpt{x}{i}{}) - \bm{g}(\idxpt{x}{i}{})\|}{\|\bm{g}(\idxpt{x}{i}{})\|},
\]

\noindent which is the discrepancy between the true gradient and the estimate, normalised by the true gradient norm. Smaller $\delta_r$ denotes a more accurate estimate.
    
The function P1, being quadratic, theoretically allows CFD to produce the exact gradients because the second-order Taylor remainders from opposing auxiliary samples cancel. However, due to the floating point (FP) representation, the decimal truncation errors of opposing auxiliary samples when $|f(\bm{x})| \;\gg 0$, and also not necessarily symmetric to cancel each other out. Therefore, such residual errors remain in the CFD numerator, and are furthermore amplified by the small step sizes, propagating to the relative error metric. To minimise the effects of FP representation in the gradient estimates (especially for $|f(\bm{x})| \;\gg 0$), we present results using \verb|float128| for our benchmarking code, which nevertheless should not be seen as a guarantee that CFD attains machine-precision accuracy at every query point.

\subsubsection*{\it Results: using nominal step sizes}

\begin{figure}
    \centering
    \includegraphics[width=\textwidth]{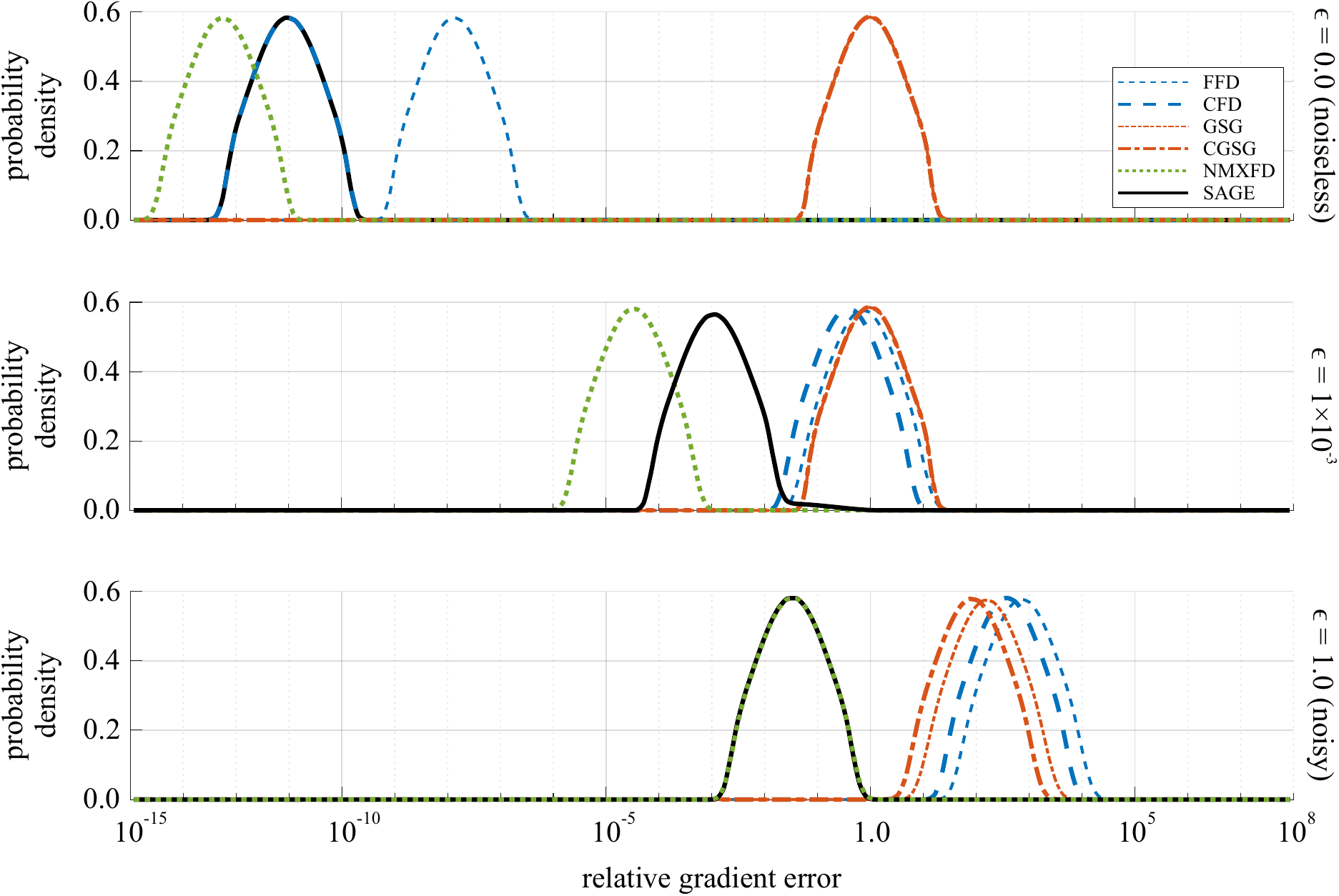}
    \caption{\textbf{P1}: relative gradient error pdf, $\kappa=1\times 10^8$, $\overline{\epsilon} = \{0.0, 1\times 10^{-3}, 1.0\}$ \\ (using nominal auxiliary step sizes)}
    \label{fig:p1-rga-fixed}
\end{figure}

\begin{figure}
    \centering
    \includegraphics[width=\textwidth]{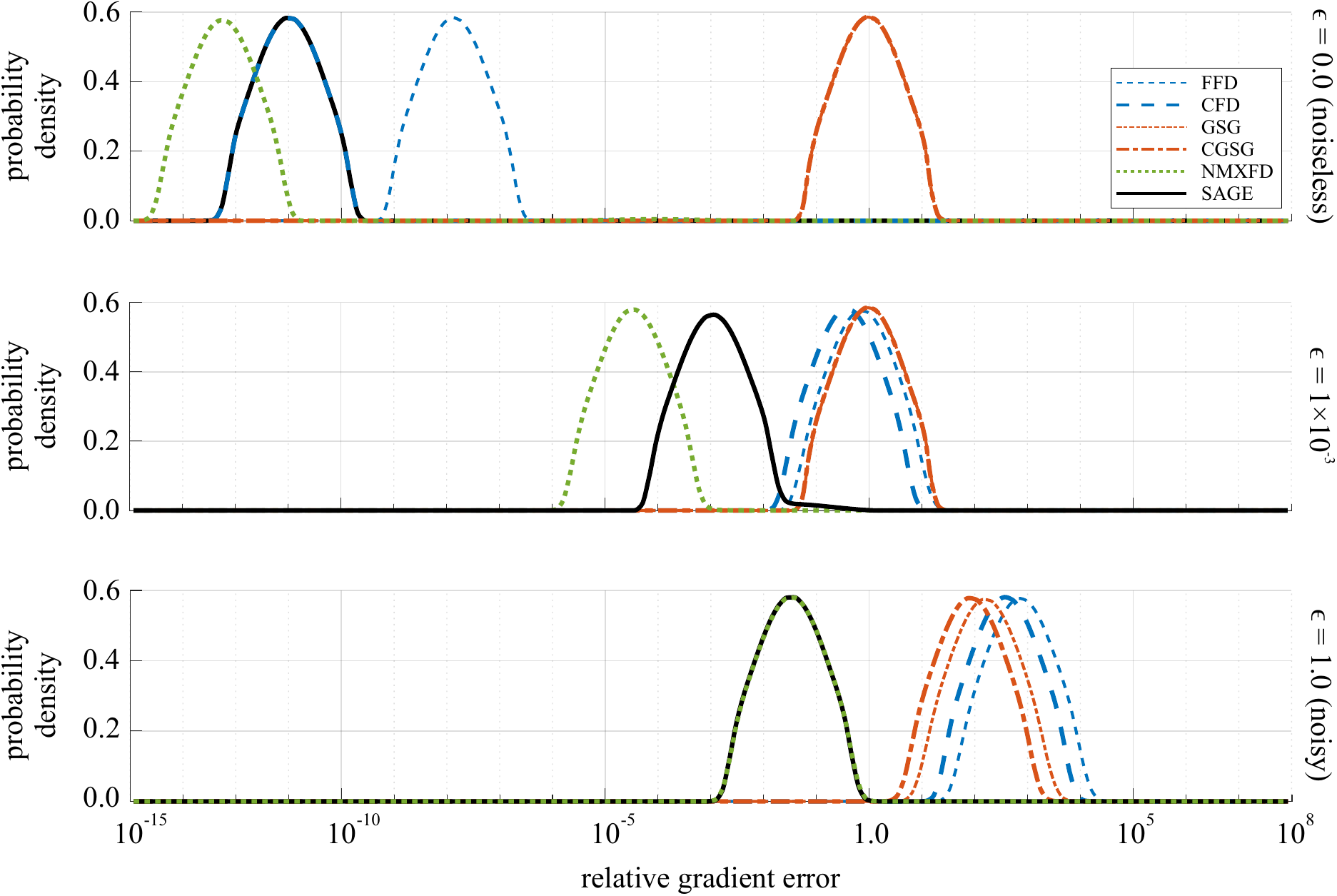}
    \caption{\textbf{P2}: relative gradient error pdf, $\kappa=1\times 10^8$, $\overline{\epsilon} = \{0.0, 1\times 10^{-3}, 1.0\}$ \\ (using nominal auxiliary step sizes)}
    \label{fig:p2-rga-fixed}
\end{figure}

\begin{figure}
    \centering
    \includegraphics[width=\textwidth]{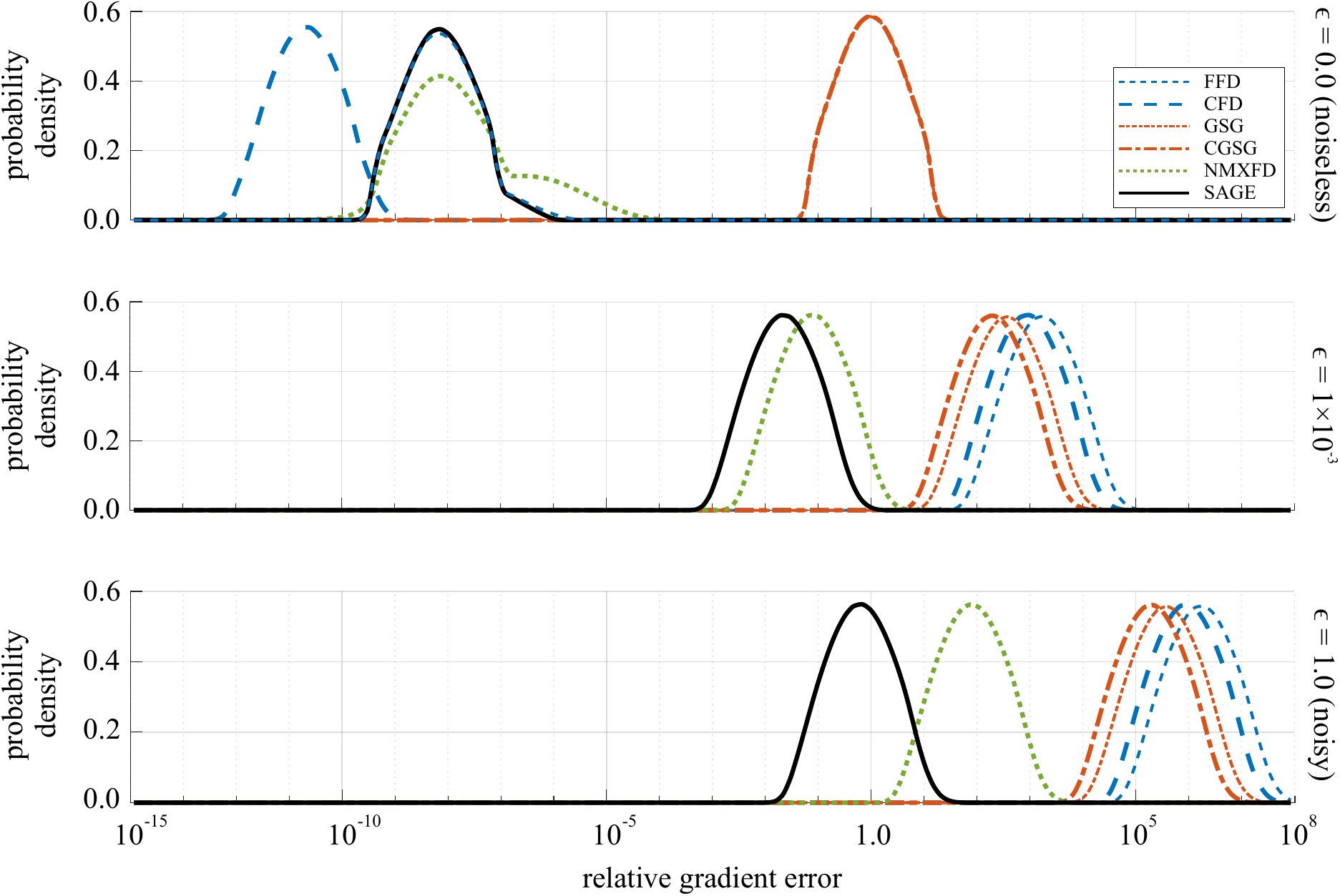}
    \caption{\textbf{P3}: relative gradient error pdf, $\kappa=1\times 10^8$, $\overline{\epsilon} = \{0.0, 1\times 10^{-3}, 1.0\}$ \\ (using nominal auxiliary step sizes)}
    \label{fig:p3-rga-fixed}
\end{figure}

\begin{figure}
    \centering
    \includegraphics[width=\textwidth]{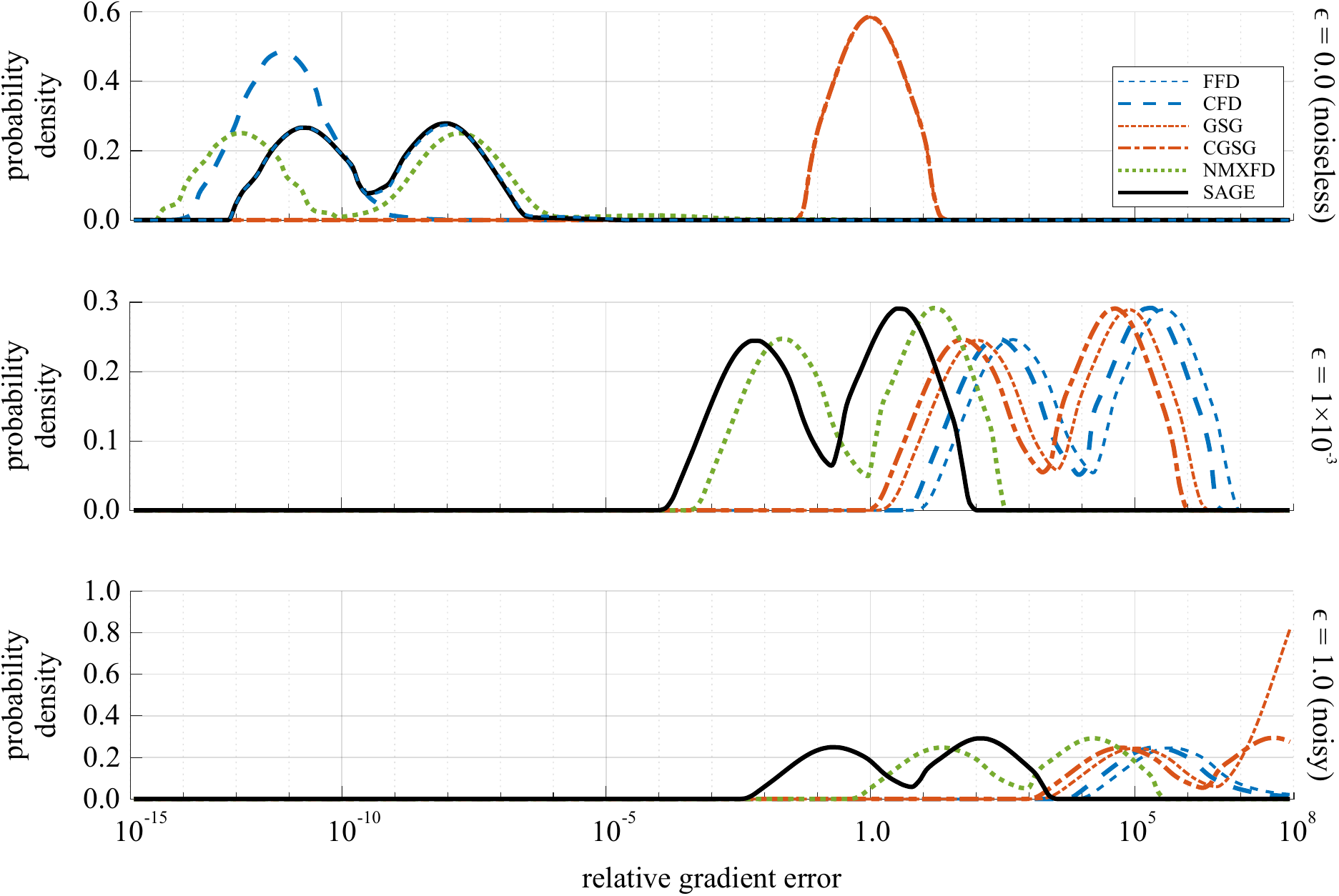}
    \caption{\textbf{P4}: relative gradient error pdf, $\kappa=1\times 10^8$, $\overline{\epsilon} = \{0.0, 1\times 10^{-3}, 1.0\}$ \\ (using nominal auxiliary step sizes)}
    \label{fig:p4-rga-fixed}
\end{figure}

\begin{figure}
    \centering
    \includegraphics[width=\textwidth]{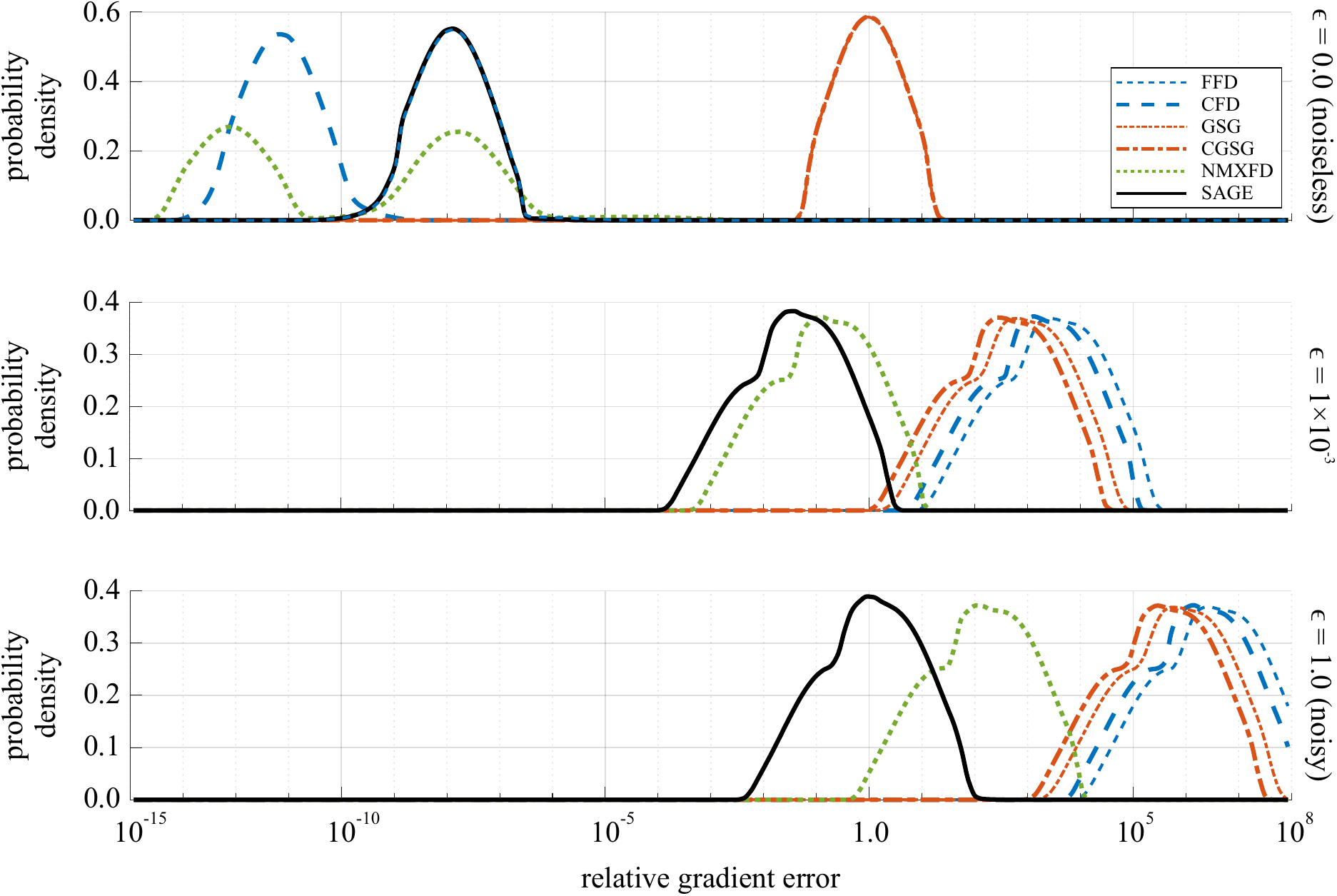}
    \caption{\textbf{P5}: relative gradient error pdf, $\kappa=1\times 10^8$, $\overline{\epsilon} = \{0.0, 1\times 10^{-3}, 1.0\}$ \\ (using nominal auxiliary step sizes)}
    \label{fig:p5-rga-fixed}
\end{figure}

We present the gradient accuracy results for the gradient estimators using their nominal step sizes (see Table~\ref{table:benchmark-gradient-estimators}). This means that the noise bounds $\overline{\epsilon}$ are used solely for generating the function values: they are not \textit{a priori} known, nor supplied as a parameter to the gradient estimators. Figs.~\ref{fig:p1-rga-fixed}-\ref{fig:p5-rga-fixed} show the probability density functions (pdf) for $\delta_r$ of the compared estimators on P1-P5, with their trends on different $\overline{\epsilon}$. In the noiseless case ($\overline{\epsilon}=0$), GSG and CGSG were consistently the least accurate estimators across all five problems, overlapping with each other and saturating around $\delta_r=1.0$ regardless of problem structure.

Among the remaining four estimators, their performances depend on the problem. On P1 and P2, CFD attained a pdf around $\delta_r \approx 1\times 10^{-11}$, consistent with theory for quadratic functions: the second-order Taylor remainders from opposite auxiliary samples cancel for CFD, and the residual $1\times 10^{-11}$ error is attributable to FP-based truncation (as discussed above). Whilst SAGE mostly coincided with CFD in these problems, NMXFD attained accuracies an order of magnitude better, at the cost of additional samples for multi-scale smoothing and a much larger nominal step size. FFD, lacking CFD's cancellation advantage on these (near-)quadratic problems, was markedly less accurate, around $\delta_r = 1\times 10^{-8}$. On P3, P4, and P5, the roles of FFD and CFD are reversed: FFD instead coincided with SAGE, whilst CFD stood alone as the most accurate estimator.

Increasing the noise bound moved the pdfs to the right for every estimator, which is unsurprising. On P1-P2, this erased CFD's noiseless cancellation advantage: once $\overline{\epsilon}>0$, the noise-driven error term dominated the now-irrelevant Taylor remainder, and CFD separated from SAGE/NMXFD, converging towards FFD's accuracy instead. NMXFD also degraded substantially with increasing noise, whilst GSG and CGSG degraded comparatively less in relative terms. SAGE, in contrast, retained the most robust accuracy as $\overline{\epsilon}$ increases. From Figs.~\ref{fig:p1-rga-fixed}-\ref{fig:p2-rga-fixed}, SAGE and NMXFD converged to essentially identical pdfs on P1 and P2 in the noisy case ($\overline{\epsilon}=1.0$). Furthermore, SAGE remained competitive or best on the other three problems, despite its initial auxiliary step size of $1\times 10^{-6}$ and lack of \textit{a priori} noise information. This highlights the necessity to adapt the sampling radius when dealing with noisy samples, which SAGE performs automatically as new samples are introduced.

Table~\ref{table:num-evals-sage} shows the number of evaluations for SAGE to arrive at the gradient queries. In the noiseless cases ($\overline{\epsilon}=0$), SAGE directly used the initial bootstrap ($2D+1$ samples) to compute the gradient estimates. However, in the noisy cases, and without noise information, the initial (uninformative) auxiliary samples (especially with initial step sizes of $1\times 10^{-6}$) are supplemented by more informative samples in the appropriate sampling radius. This raised the evaluation counts to 61-81 (depending on the problem and noise bounds), which is an expected cost of estimating $\overline{\epsilon}$ directly from data. However, as demonstrated in the next benchmarks, such a need for additional samples is compensated by samples reuse, leading to better optimisation iterate improvements with limited evaluations budget.

\begin{table}[]
\small
\centering
\begin{tabular}{l|cc|cc|cc|cc|cc}
\toprule
\textbf{Noise $\overline{\epsilon}$}    & \multicolumn{2}{c}{\textbf{P1}} & \multicolumn{2}{c}{\textbf{P2}} & \multicolumn{2}{c}{\textbf{P3}} & \multicolumn{2}{c}{\textbf{P4}} & \multicolumn{2}{c}{\textbf{P5}} \\
\midrule
0.0 & 41 & 41       & 41 & 41       & 41 & 41       & 41 & 41       & 41 & 41 \\
$1\times 10^{-3}$ & 61 & 61    & 61 & 61       & 61 & 61       & 62 & 81       & 62 & 81   \\
1.0 & 61 & 61       & 61 & 61       & 81 & 81      & 62 & 81      & 62 & 81  \\
\bottomrule
\end{tabular}
\caption{SAGE number of evaluations for each gradient query \\ (median / p90, lower is better)}
\label{table:num-evals-sage}
\end{table}

\subsubsection*{\it Results: using noise-adjusted auxiliary step sizes}

We additionally considered the case with \textit{a priori} known $\overline{\epsilon}$, and the auxiliary step sizes of FFD, CFD, GSG, and CGSG are adjusted accordingly (now denoted FFD-SSR, CFD-SSR, GSG-SSR, and CGSG-SSR). On the other hand, neither NMXFD nor SAGE receives a step size adjustment: no step-size rule is established for NMXFD \cite{boresta2022mixed}, and SAGE continues to estimate $\tilde{\epsilon}$ from samples exactly as in Figs.~\ref{fig:p1-rga-fixed}-\ref{fig:p5-rga-fixed}. This is consistent with our premise of SAGE operating without \textit{a priori} noise characterisation. Therefore, NMXFD and SAGE are reproduced from Figs.~\ref{fig:p1-rga-fixed}-\ref{fig:p5-rga-fixed} for reference. According to the step size rules found in \cite{berahas2022comparison}, we pre-adjusted the step sizes of FFD, CFD, GSG, and CGSG for $\overline{\epsilon}=\{1\times 10^{-3}, 1.0\}$; at $\overline{\epsilon}=0$ these rules reduce to the nominal values of Table~\ref{table:benchmark-gradient-estimators}, so this comparison adds information only at the two noisy levels.

\begin{figure}
    \centering
    \includegraphics[width=\textwidth]{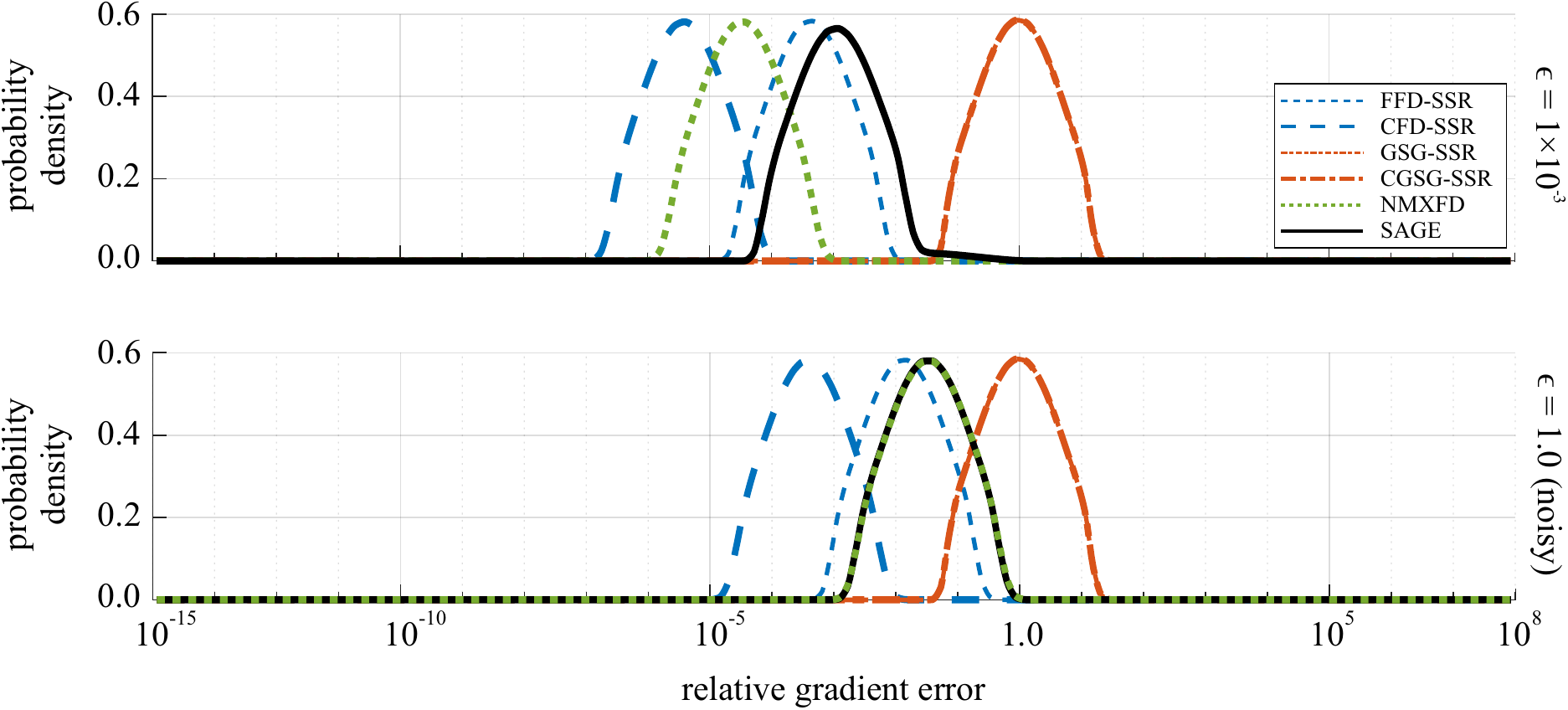}
    \caption{\textbf{P1}: relative gradient error pdf, $\kappa=1\times 10^8$, $\overline{\epsilon} = \{0.0, 1\times 10^{-3}, 1.0\}$ \\ (using noise-adjusted auxiliary step sizes)}
    \label{fig:p1-rga-adj}
\end{figure}

\begin{figure}
    \centering
    \includegraphics[width=\textwidth]{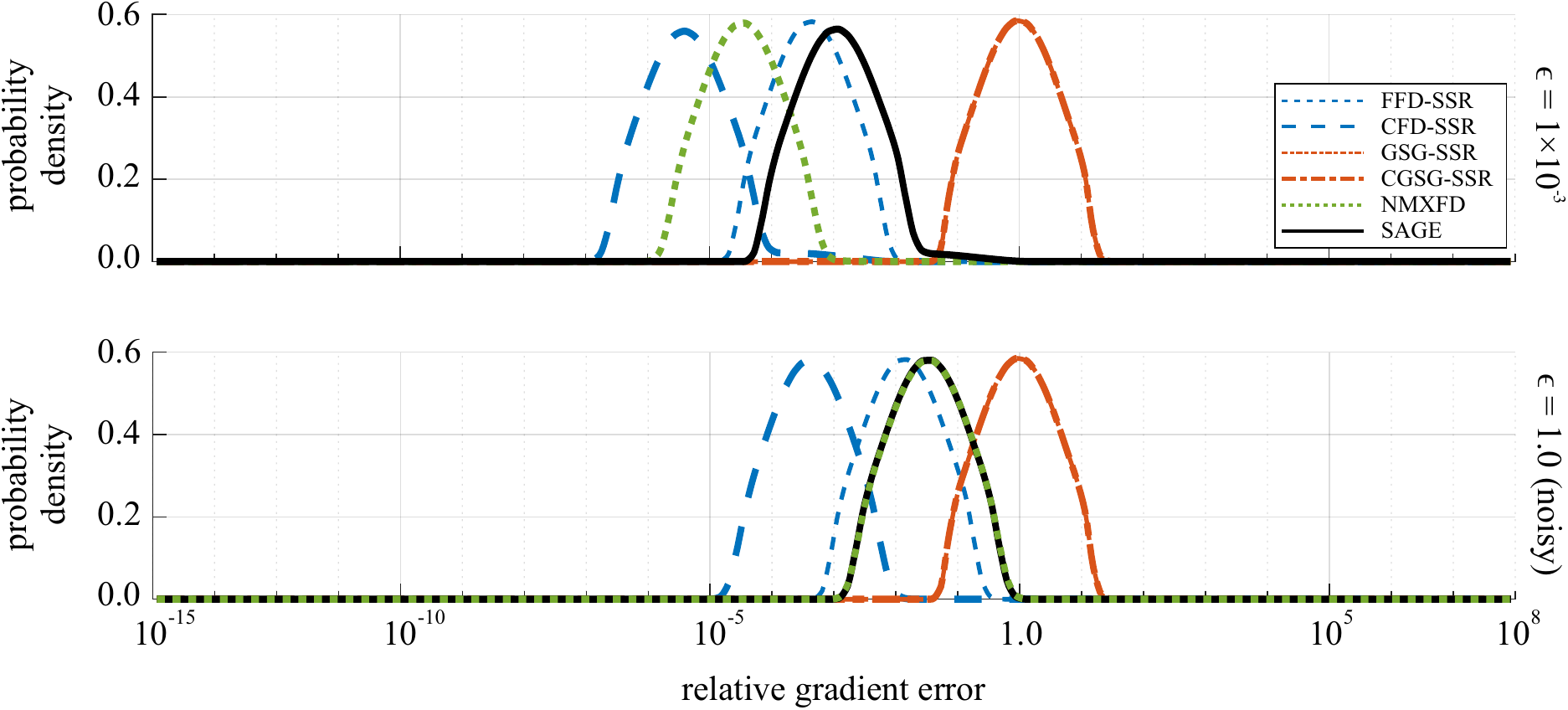}
    \caption{\textbf{P2}: relative gradient error pdf, $\kappa=1\times 10^8$, $\overline{\epsilon} = \{0.0, 1\times 10^{-3}, 1.0\}$ \\ (using noise-adjusted auxiliary step sizes)}
    \label{fig:p2-rga-adj}
\end{figure}

\begin{figure}
    \centering
    \includegraphics[width=\textwidth]{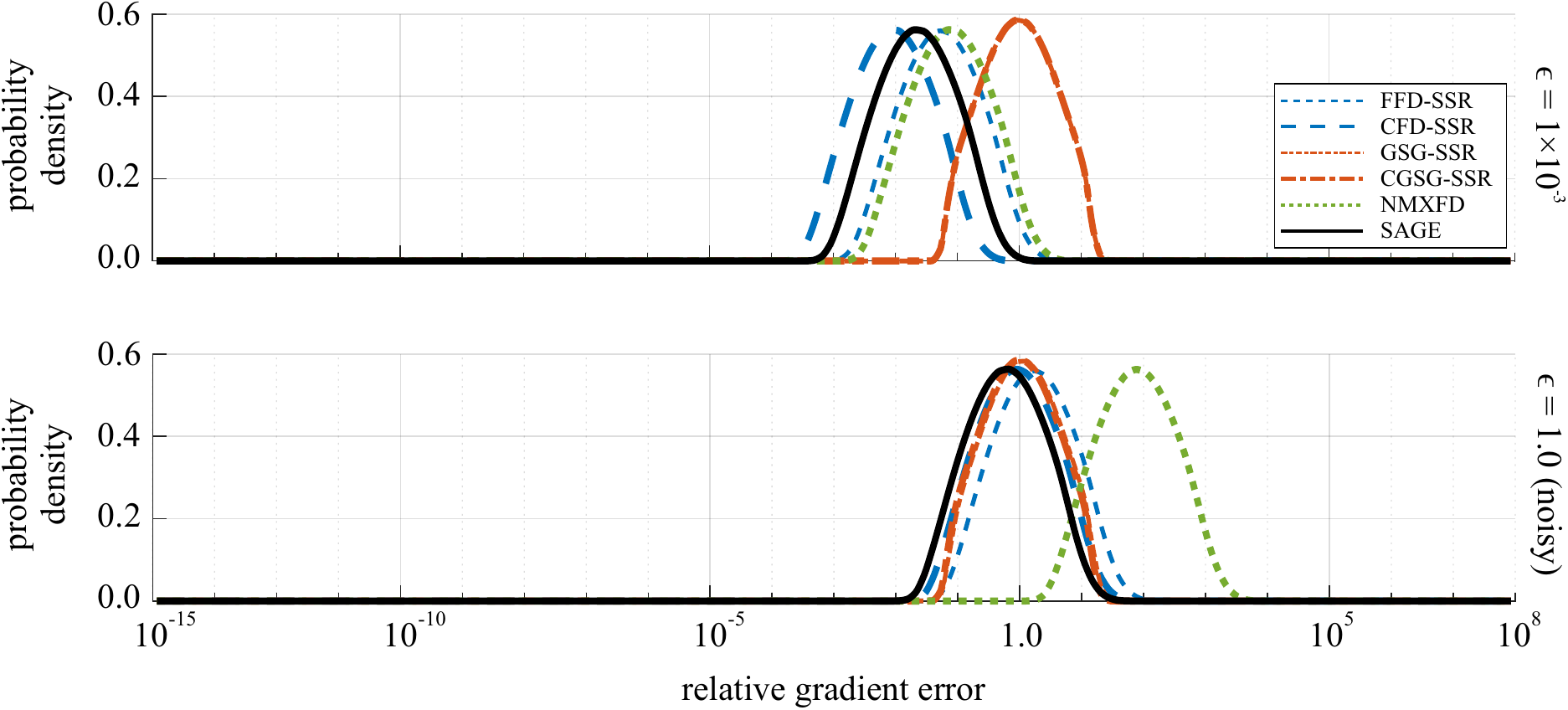}
    \caption{\textbf{P3}: relative gradient error pdf, $\kappa=1\times 10^8$, $\overline{\epsilon} = \{0.0, 1\times 10^{-3}, 1.0\}$ \\ (using noise-adjusted auxiliary step sizes)}
    \label{fig:p3-rga-adj}
\end{figure}

\begin{figure}
    \centering
    \includegraphics[width=\textwidth]{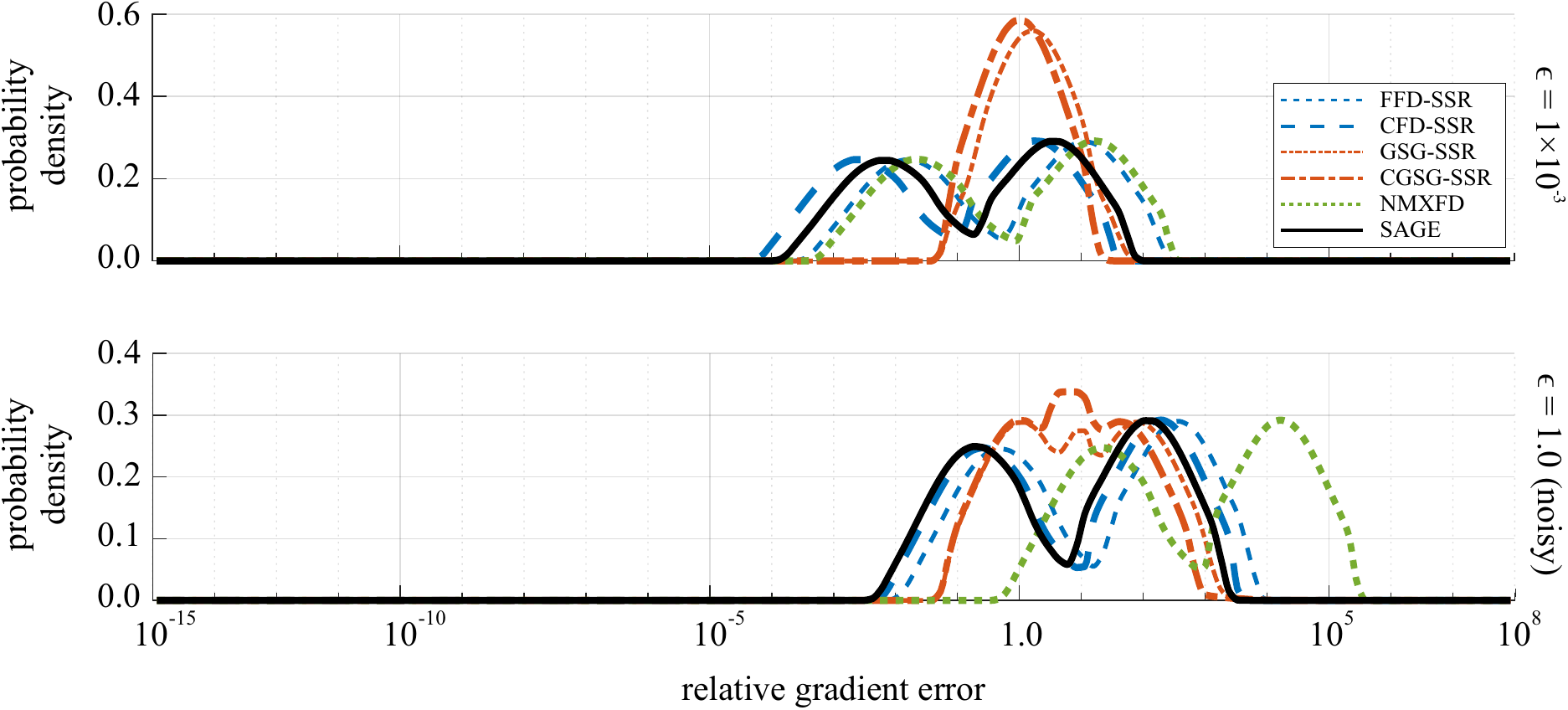}
    \caption{\textbf{P4}: relative gradient error pdf, $\kappa=1\times 10^8$, $\overline{\epsilon} = \{0.0, 1\times 10^{-3}, 1.0\}$ \\ (using noise-adjusted auxiliary step sizes)}
    \label{fig:p4-rga-adj}
\end{figure}

\begin{figure}
    \centering
    \includegraphics[width=\textwidth]{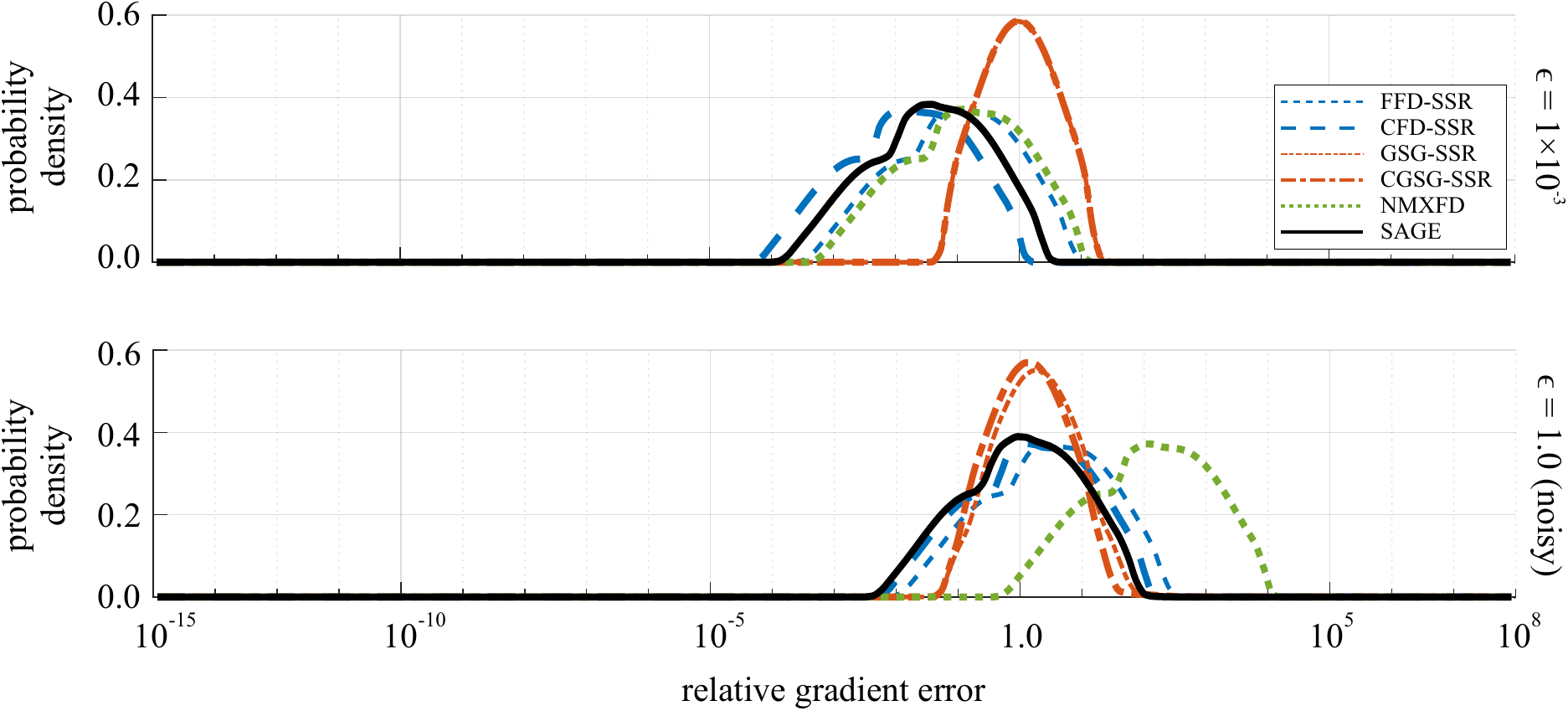}
    \caption{\textbf{P5}: relative gradient error pdf, $\kappa=1\times 10^8$, $\overline{\epsilon} = \{0.0, 1\times 10^{-3}, 1.0\}$ \\ (using noise-adjusted auxiliary step sizes)}
    \label{fig:p5-rga-adj}
\end{figure}

Figs.~\ref{fig:p1-rga-adj}-\ref{fig:p5-rga-adj} show the corresponding pdfs for P1-P5 at $\overline{\epsilon}=\{1\times 10^{-3}, 1.0\}$. As expected, the pdfs of FFD-SSR, CFD-SSR, GSG-SSR, and CGSG-SSR became more robust because of the noise-informed step sizes: FFD-SSR, in particular, now overlapped with SAGE on P1-P2, with CFD-SSR resulting in the best pdfs. Even so, on P3-P5, SAGE has either achieved competitive or better pdfs compared to all four noise-informed estimators, despite lacking the noise information they were given.

\subsection{\it First-order optimisation benchmark tests}

We tested the compared gradient estimators within a simple gradient descent-based optimisation routine. For P1-P5, we tested with a high condition number $\kappa$ for $\bm{Q}$ ($\kappa = 1\times10^8$), representing a difficult optimisation landscape. Trying different noise bounds $\overline{\epsilon}=\left\{0.0, 1\times10^{-3}, 1.0\right\}$ as in the previous benchmark, we performed 100 trials, each one using a different randomly-generated $\bm{Q}$, $\bm{y}$, and $\idxpt{x}{1}{}$, and with a budget of $N=50D$ (1000 in this paper) function evaluations. For fairness across competitor methods, all random generations were done with the same set of seeds.

All comparative tests in this benchmark used the same first-order optimisation routine, differing only in the gradient estimation module, and whose parameters are summarised in Table~\ref{table:first-order-solver-params}. Specifically, we employed a simple gradient descent method with backtracking line search, implemented once in a common solver class. The solver performs the update rule of the form

    \[
        \idxpt{x}{k+1}{} = \idxpt{x}{k}{} - \eta \tilde{\bm{g}}^{(k)},
    \]

\noindent where $\tilde{\bm{g}}^{(k)}$ is the estimated gradient at the current iterate $\idxpt{x}{k}{}$ (using the estimator under test), and $\idxval{\eta}{k}{}$ is selected via Armijo backtracking line search. Backtracking begins with an initial step size $\eta_0 = 1.0$ and reduces it by 0.5, until the Armijo condition

    \[
        f(\idxpt{x}{k}{} - \eta \tilde{\bm{g}}^{(k)} ) \leq f(\idxpt{x}{k}{}) - 1\times10^{-6} \eta \|\tilde{\bm{g}}^{(k)}\|^2
    \]

\noindent is satisfied. When $\eta \leq 1\times10^{-6}$, the line search is restarted: the gradient estimate $\tilde{\bm{g}}^{(k)}$ is recomputed (which, especially with noisy evaluations, avoids repeated line searches along stale descent directions), and $\eta$ is reset to $\eta_0$, with all evaluations performed counted towards the evaluation budget. All samples, including rejected trial points during line search, are passed to the gradient estimator. 

Optimisation runs were terminated solely on the evaluation budget $N$, counting all function evaluations (including auxiliary samples for gradient estimation, and those generated by backtracking line searches). This means that no optimality-based or progress-based stopping criteria were employed. Such a fixed-budget protocol ensures no interaction between the first-order optimisation method and the tested gradient estimator, leading to a fair comparison of descent effectiveness given the same amount of information.

\begin{table}[]
\small
\centering
\begin{tabular}{ll}
    \toprule
    \textbf{Parameter}        & \textbf{Value}                                     \\ 
    \midrule
    Update rule               & $\idxpt{x}{k+1}{} = \idxpt{x}{k}{} - \eta \tilde{\bm{g}}^{(k)}$ \\
    Step-size mode            & Adaptive (Armijo backtracking) \\
    Initial step size         & 1.0                                                \\ 
    Backtracking factor       & 0.5                                                \\ 
    Armijo parameter          & $1\times10^{-6}$ \\ 
    Minimum step size         & $1\times10^{-6}$ (reset to $\eta_0$) \\
    Stop criterion            & $n = N$ (evaluation budget exhausted) \\
    \bottomrule
\end{tabular}
\caption{First-order solver parameters}
\label{table:first-order-solver-params}
\end{table}

\subsubsection*{\it Comparison metric}

We used \textit{data profiles} as our comparison metric, following the benchmarking framework introduced in \cite{more2009benchmarking}. Data profiles quantify the probability that a method achieves a prescribed level of progress within a given number of function evaluations, aggregated across a set of problem instances.

For a given problem $p$ (P1-P5) and the first-order method summarised in Table~\ref{table:first-order-solver-params} (with gradient estimator $m$), let $n_{pm}$ denote the number of function evaluations required to achieve a target relative improvement (and thereby marked as ``solved''). In our paper, we set this goal as $\frac{f(\idxpt{x}{n}{})}{f(\idxpt{x}{1}{})} \leq 0.01$, and if this is not reached within the maximum evaluations $N$, we set $n_{pm} = +\infty$. The data profile after $n$ evaluations is then defined as

\[
    d_m(\tau) = \frac{1}{|\mathcal{P}|} \left| \left\{ p \in \mathcal{P} \; : \; n_{pm} \le \tau \right\} \right|,
\]

\noindent where $\mathcal{P}$ is the set of all tested problem instances and $\tau = \frac{n}{N}$ is a normalised evaluation budget. Thus, $d_m(\tau)$ is the fraction of problems solved within budget $\tau$. In tests with noisy samples, the relative improvement is evaluated using the corresponding \textit{noise-free} objective values at the sampled iterates, which are available for the benchmark problems (i.e., not to the gradient estimators) and used only for post-hoc analysis.

\subsubsection*{Results: effect of noise levels (using nominal step sizes)}

\begin{figure}
    \includegraphics[width=\textwidth]{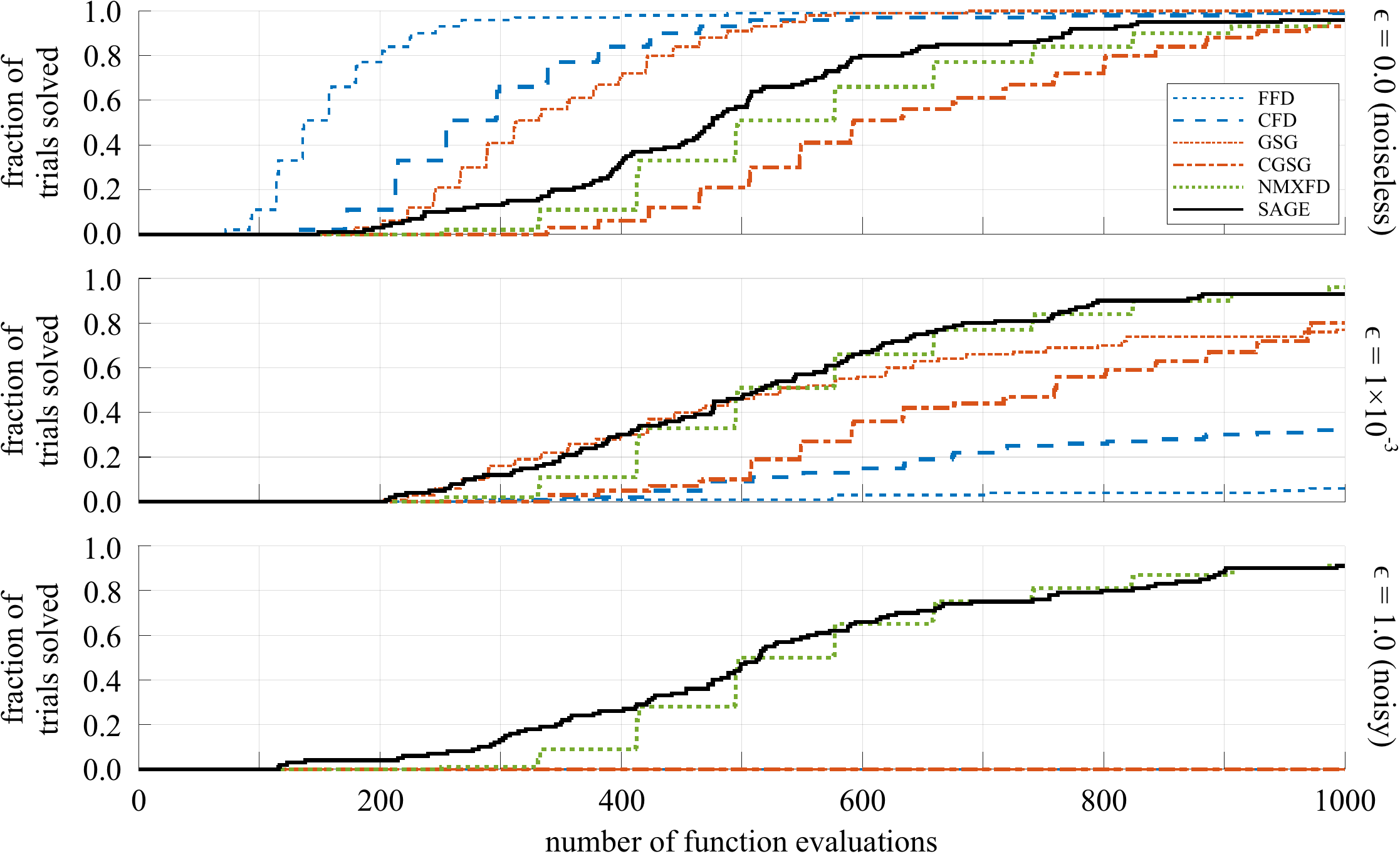}
    \caption{\textbf{P1}: data profiles, $\kappa=1\times 10^8$, $\overline{\epsilon} = \{0.0, 1\times 10^{-3}, 1.0\}$ \\ (using nominal auxiliary step sizes)}
    \label{fig:data-profiles-lsq-varying-noise}
\end{figure}

\begin{figure}
    \includegraphics[width=\textwidth]{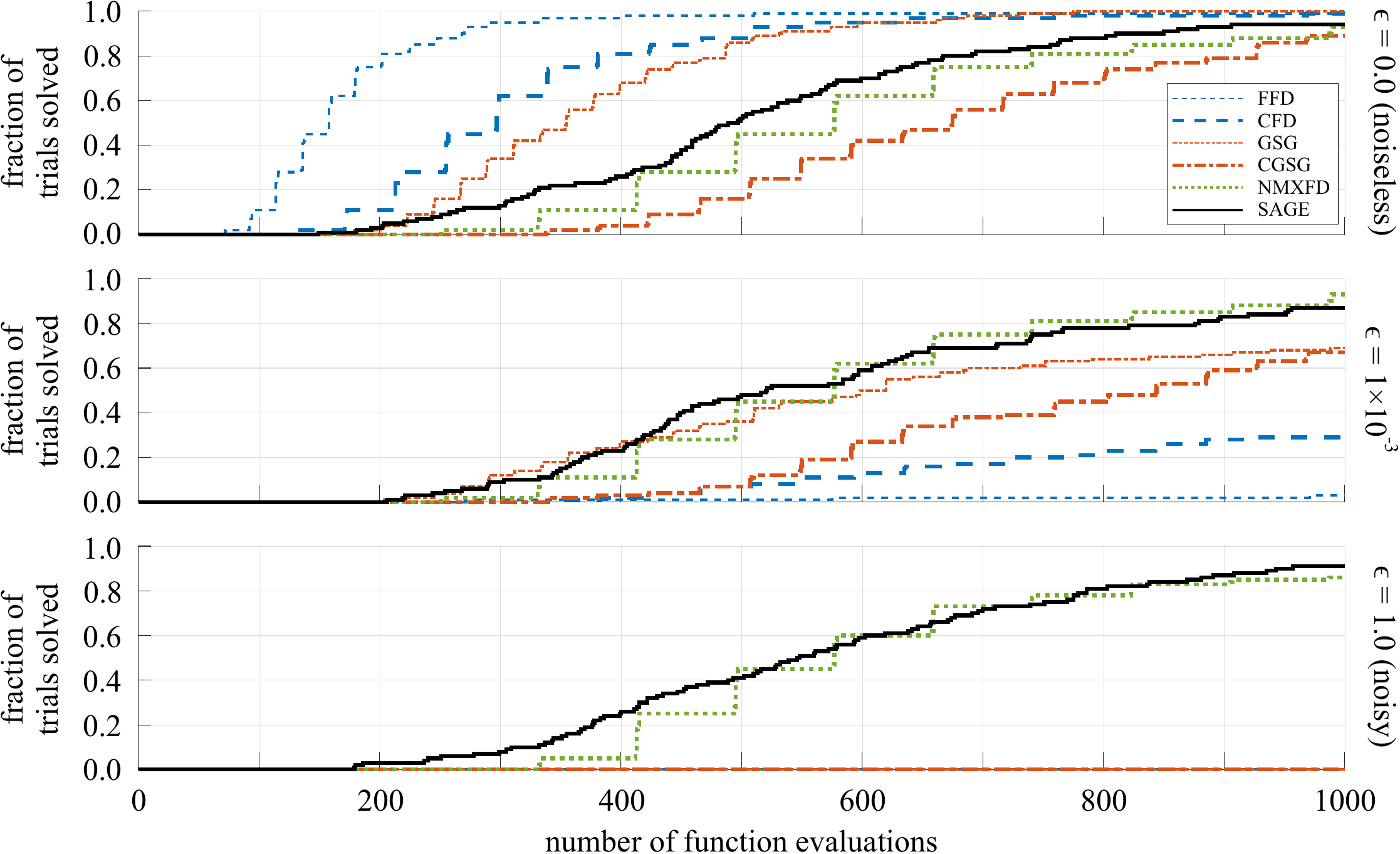}
    \caption{\textbf{P2}: data profiles, $\kappa=1\times 10^8$, $\overline{\epsilon} = \{0.0, 1\times 10^{-3}, 1.0\}$ \\ (using nominal auxiliary step sizes)}
    \label{fig:data-profiles-lasso-varying-noise}
\end{figure}

\begin{figure}
    \includegraphics[width=\textwidth]{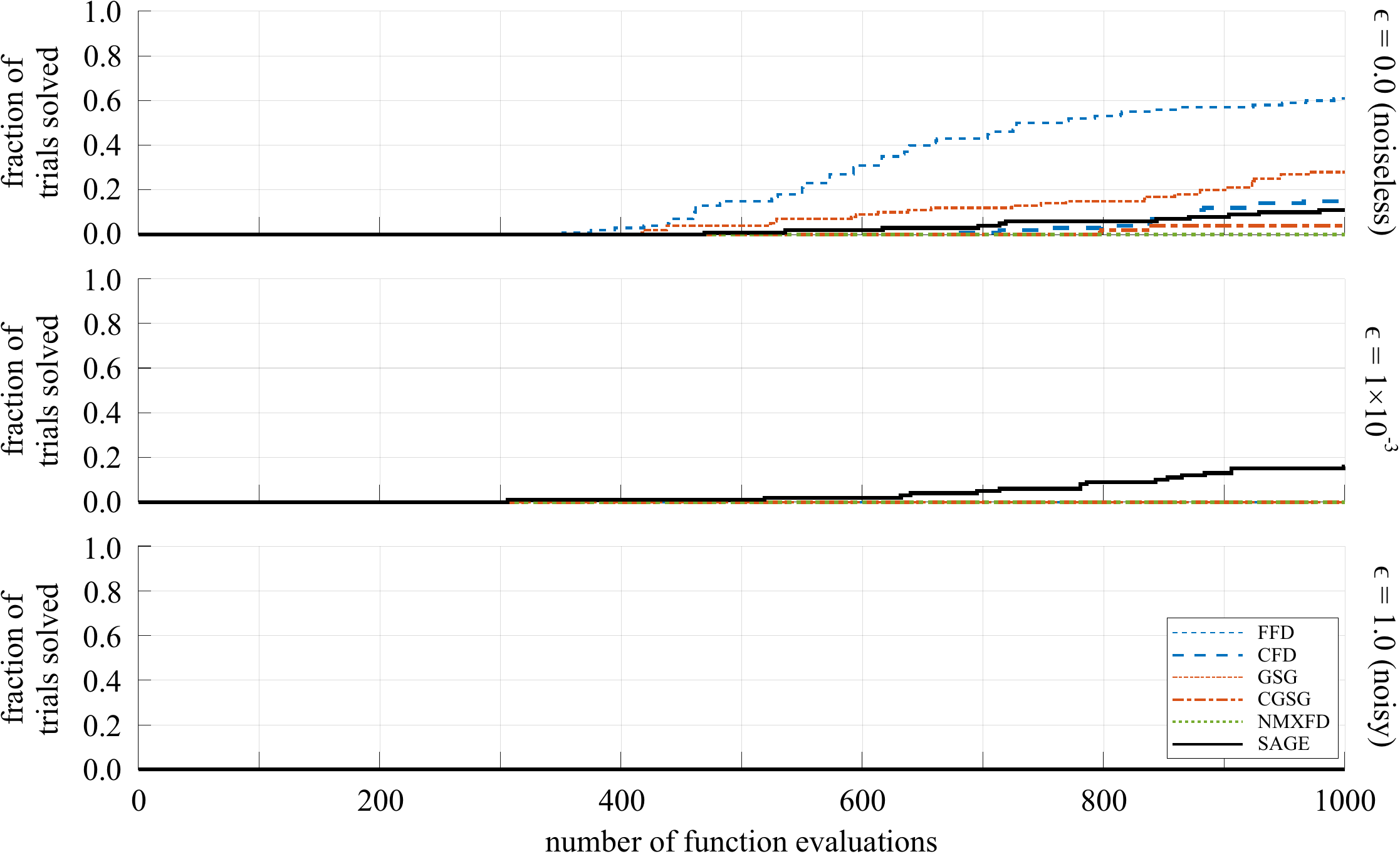}
    \caption{\textbf{P3}: data profiles, $\kappa=1\times 10^8$, $\overline{\epsilon} = \{0.0, 1\times 10^{-3}, 1.0\}$ \\ (using nominal auxiliary step sizes)}
    \label{fig:data-profiles-logsumexp-varying-noise}
\end{figure}

\begin{figure}
    \includegraphics[width=\textwidth]{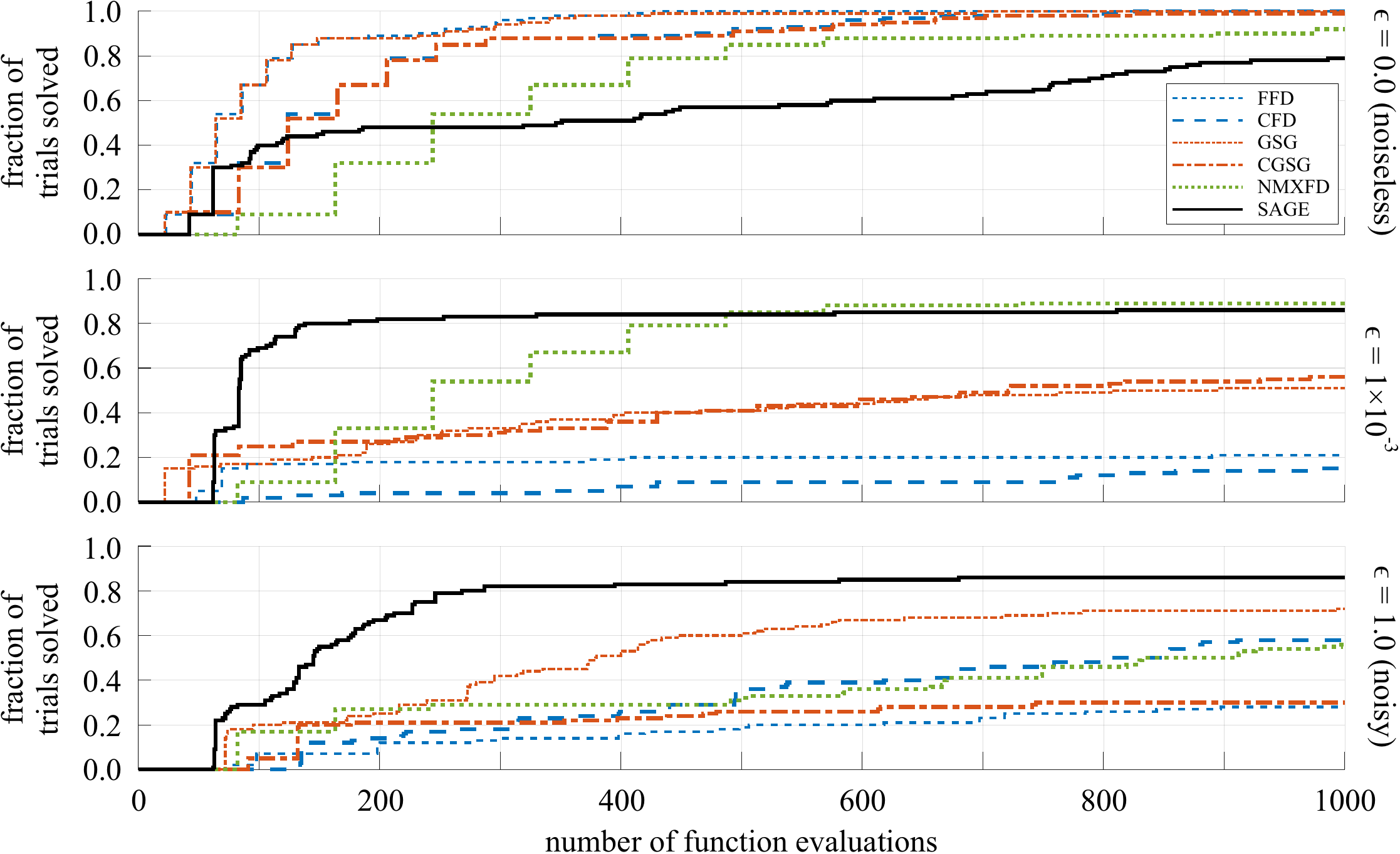}
    \caption{\textbf{P4}: data profiles, $\kappa=1\times 10^8$, $\overline{\epsilon} = \{0.0, 1\times 10^{-3}, 1.0\}$ \\ (using nominal auxiliary step sizes)}
    \label{fig:data-profiles-l1logreg-varying-noise}
\end{figure}

\begin{figure}
    \includegraphics[width=\textwidth]{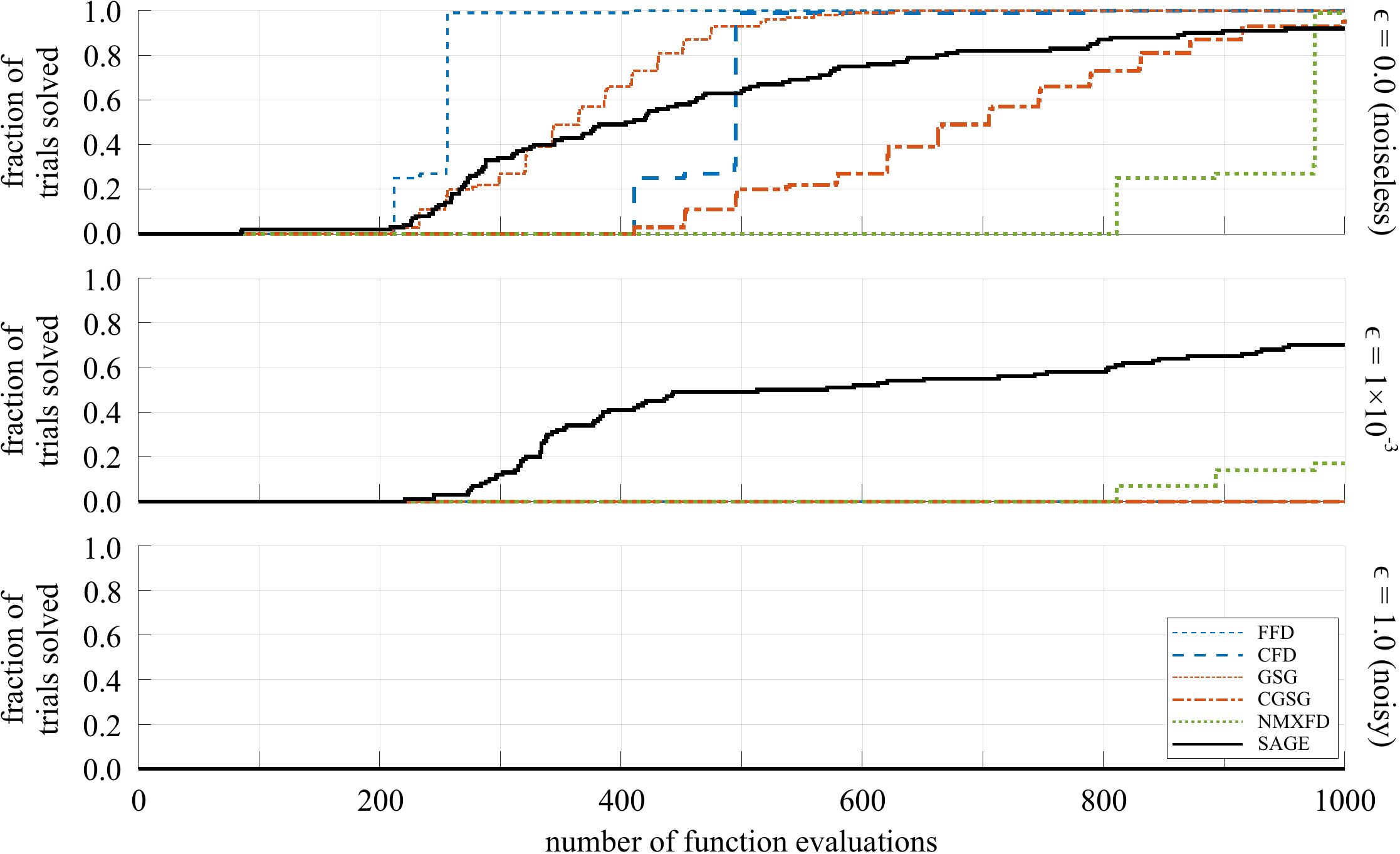}
    \caption{\textbf{P5}: data profiles, $\kappa=1\times 10^8$, $\overline{\epsilon} = \{0.0, 1\times 10^{-3}, 1.0\}$ \\ (using nominal auxiliary step sizes)}
    \label{fig:data-profiles-l2logreg-varying-noise}
\end{figure}

Figs.~\ref{fig:data-profiles-lsq-varying-noise}-\ref{fig:data-profiles-l2logreg-varying-noise} show that as noise bound $\overline{\epsilon}$ increased, the optimisation performance degraded among the compared methods, with worsened data profiles as expected. In most cases, the data profiles shifted rightward and the final success rates decreased with higher noise (particularly at $\overline{\epsilon}=1.0$). For P1 and P2, the data profiles along with their relative order in performance are very similar, with only SAGE and NMXFD achieving ``solved'' runs (i.e., attaining the desired improvement), with SAGE achieving better data profiles. This can be traced to the much higher evaluation requirements per iteration for NMXFD, such that only a few iterations were done within the evaluation budget. In P4, SAGE fared the worst in the noiseless case, but emerged with the best data profiles under both noisy settings. On the other hand, for P3 and P5, even modest noise ($\overline{\epsilon}=1\times10^{-3}$) substantially degraded the data profiles of FFD, CFD, GSG, CGSG, and NMXFD, with only SAGE remaining with solved runs on P3; no method had any solved runs on both problems with $\overline{\epsilon}=1.0$. We also tested with (unbounded) Gaussian noise, and the results were consistent with what we show here with bounded noise.

The performance degradation with increasing $\overline{\epsilon}$, especially for finite-difference-based estimators that rely on fixed perturbation scales or smoothing parameters, is unsurprising. As the noise bound increases, small finite-difference steps amplified the noise introduced, leading to less reliable gradient information under a fixed evaluation budget. This behaviour is reflected by the rapid drop in success rates for FFD and CFD as noise increases, especially on problems P3 and P5. Smoothing-based estimators such as GSG and CGSG exhibited a similar sensitivity, with performance deteriorating as noise overwhelms the averaging effect of a limited number of samples. In contrast, NMXFD and SAGE retained solved runs at higher noise levels on several problems, reflecting their higher per-iteration evaluation cost, and more conservative use of samples information, respectively.

\subsubsection*{Results: effect of noise levels (using noise-adjusted step sizes)}

\begin{figure}
    \includegraphics[width=\textwidth]{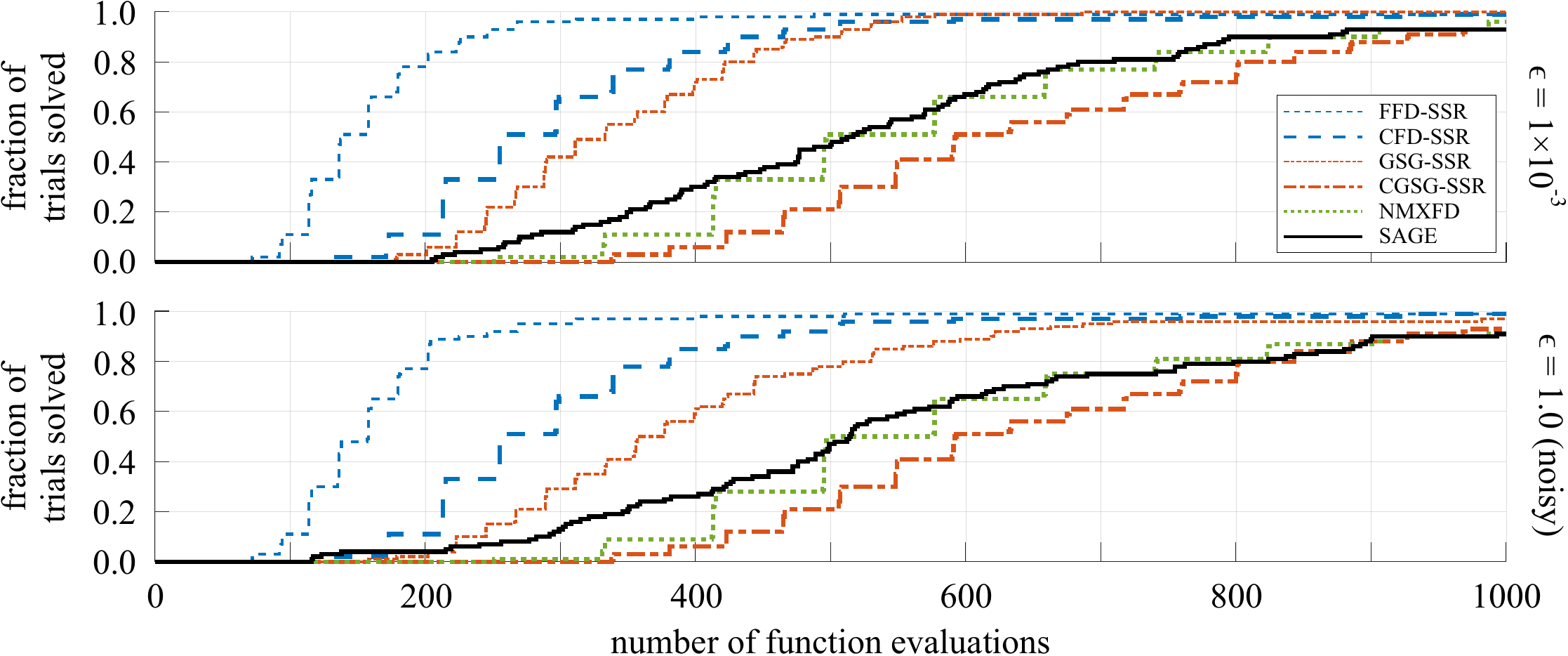}
    \caption{\textbf{P1}: data profiles, $\kappa=1\times 10^8$, $\overline{\epsilon} = \{0.0, 1\times 10^{-3}, 1.0\}$ \\ (using noise-adjusted auxiliary step sizes)}
    \label{fig:data-profiles-lsq-varying-noise2}
\end{figure}

\begin{figure}
    \includegraphics[width=\textwidth]{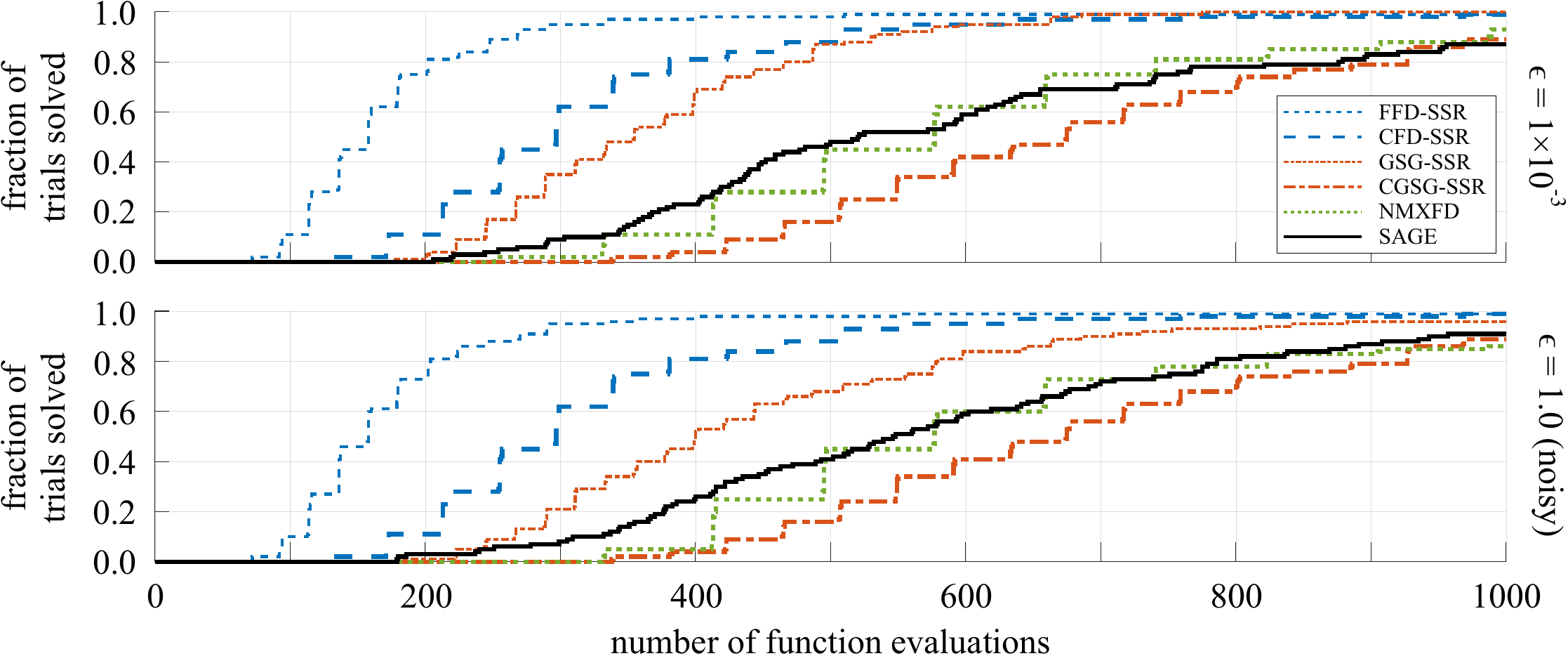}
    \caption{\textbf{P2}: data profiles, $\kappa=1\times 10^8$, $\overline{\epsilon} = \{0.0, 1\times 10^{-3}, 1.0\}$ \\ (using noise-adjusted auxiliary step sizes)}
    \label{fig:data-profiles-lasso-varying-noise2}
\end{figure}

\begin{figure}
    \includegraphics[width=\textwidth]{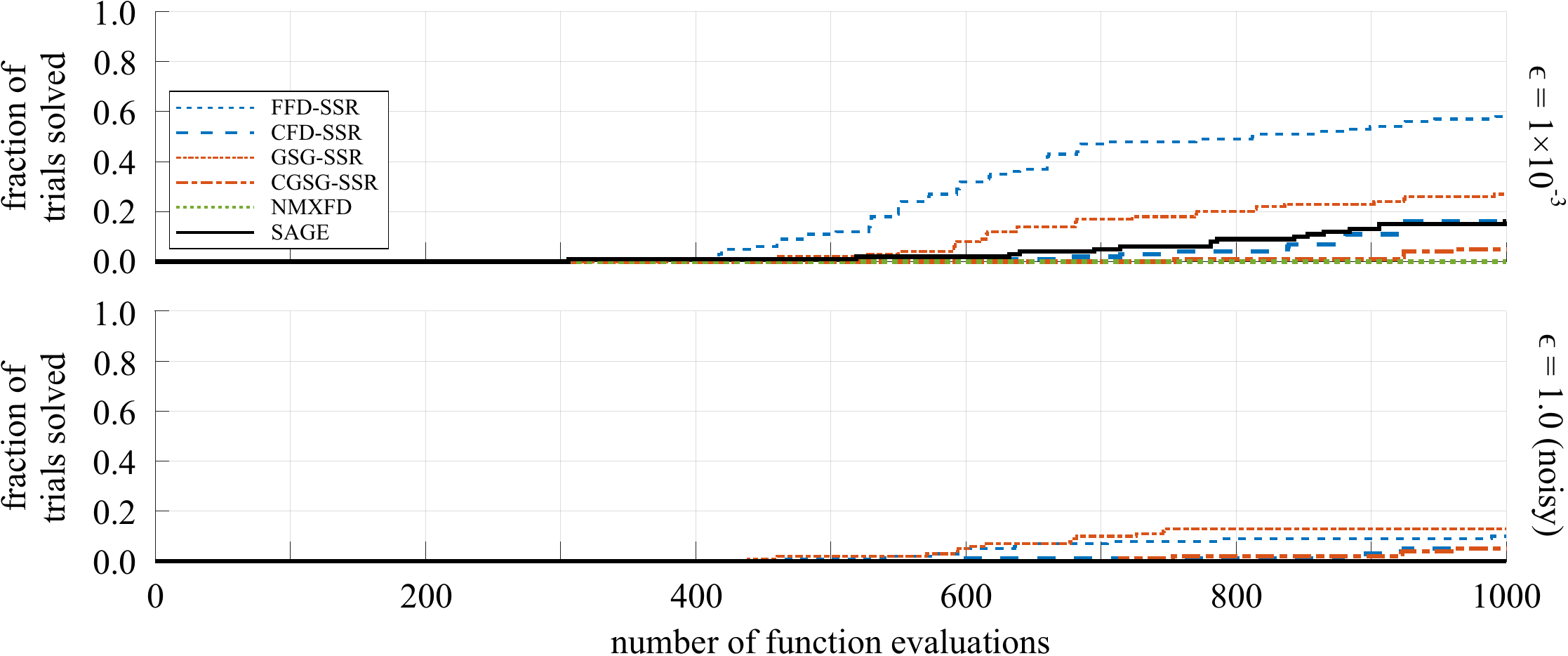}
    \caption{\textbf{P3}: data profiles, $\kappa=1\times 10^8$, $\overline{\epsilon} = \{0.0, 1\times 10^{-3}, 1.0\}$ \\ (using noise-adjusted auxiliary step sizes)}
    \label{fig:data-profiles-logsumexp-varying-noise2}
\end{figure}

\begin{figure}
    \includegraphics[width=\textwidth]{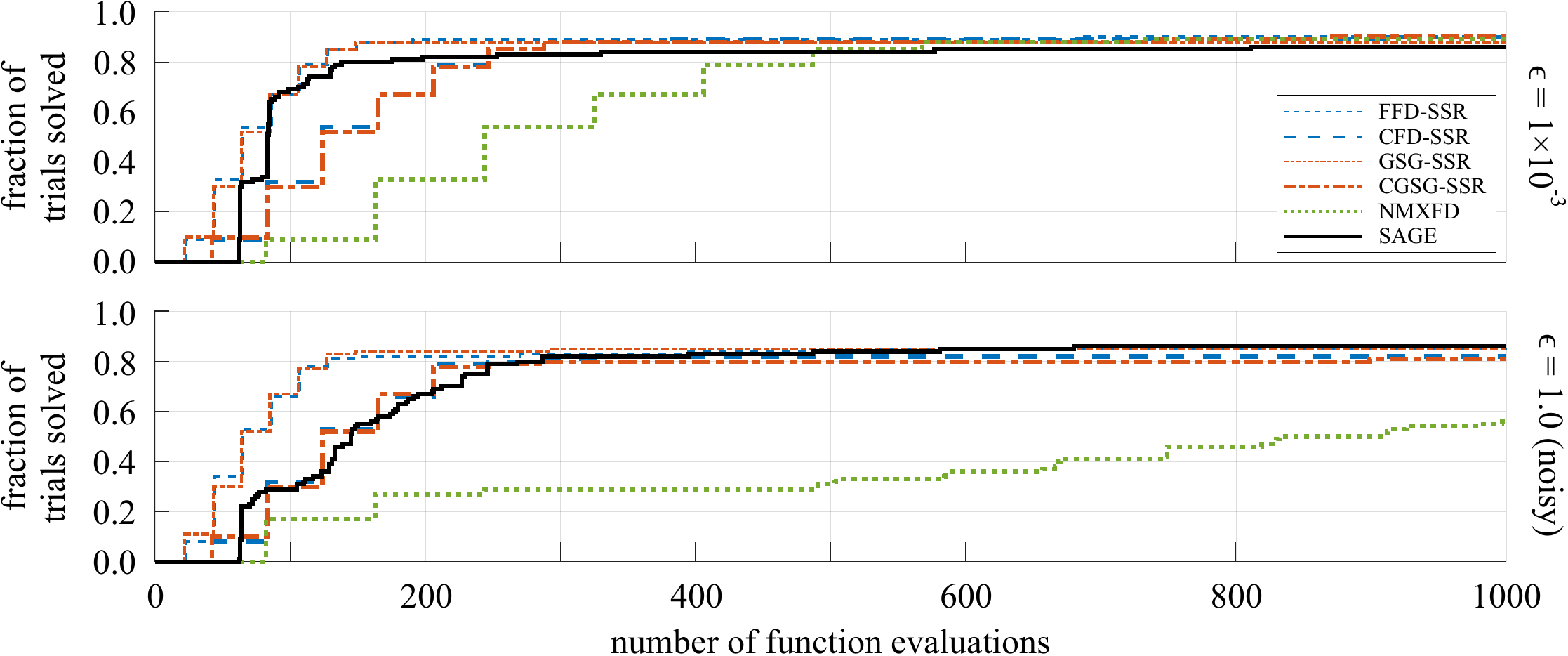}
    \caption{\textbf{P4}: data profiles, $\kappa=1\times 10^8$, $\overline{\epsilon} = \{0.0, 1\times 10^{-3}, 1.0\}$ \\ (using noise-adjusted auxiliary step sizes)}
    \label{fig:data-profiles-l1logreg-varying-noise2}
\end{figure}

\begin{figure}
    \includegraphics[width=\textwidth]{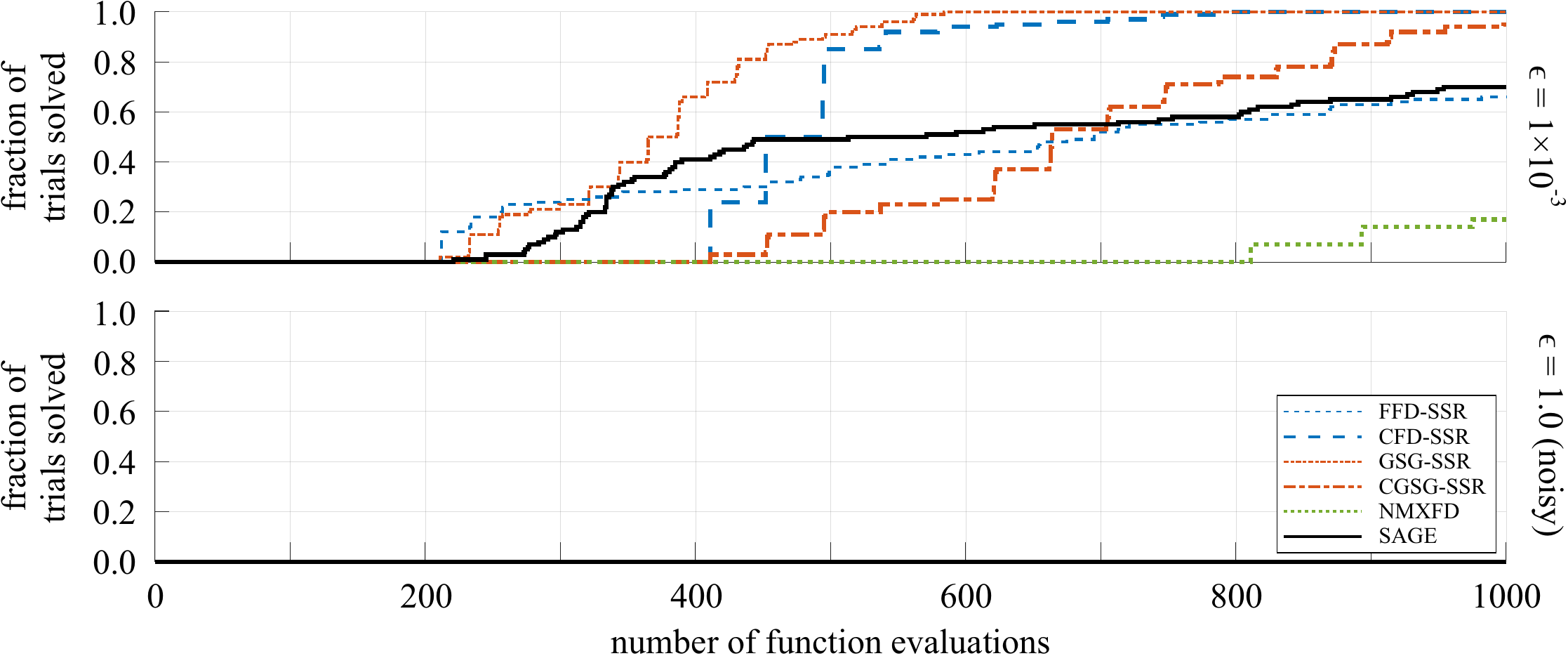}
    \caption{\textbf{P5}: data profiles, $\kappa=1\times 10^8$, $\overline{\epsilon} = \{0.0, 1\times 10^{-3}, 1.0\}$ \\ (using noise-adjusted auxiliary step sizes)}
    \label{fig:data-profiles-l2logreg-varying-noise2}
\end{figure}

We repeated the first-order optimisation-based comparison, but now with FFD, CFD, GSG, and CGSG using the noise-informed step size rules (denoted again as FFD-SSR, CFD-SSR, GSG-SSR, and CGSG-SSR), as a sanity check similar to the previous gradient accuracy benchmarks. As before, SAGE and NMXFD did not have a comparable step size adjustment. Shown in Figs.~\ref{fig:data-profiles-lsq-varying-noise2}-\ref{fig:data-profiles-l2logreg-varying-noise2} are the data profiles for P1-P5 at $\overline{\epsilon}=\{1\times 10^{-3}, 1.0\}$ (the noiseless case is omitted in these plots because they reduce anyway to the nominal step sizes, whose plots are in the top subplots of Figs.~\ref{fig:data-profiles-lsq-varying-noise}-\ref{fig:data-profiles-l2logreg-varying-noise}).

As seen in the data profiles, the outcome varies by problem. On P1 and P2, FFD-SSR reached the best data profiles at both noise levels. SAGE eventually attained a comparable final success rate as all other estimators within the evaluation budget, but is clearly less efficient than FFD-SSR, CFD-SSR, and GSG-SSR. On P4, SAGE was broadly comparable to FFD-SSR, GSG-SSR on moderate noise, and CFD-SSR and CGSG-SSR at high noise. On P5 at $\overline{\epsilon}=1\times10^{-3}$, SAGE was comparable to FFD-SSR, but was clearly outperformed by CFD-SSR, GSG-SSR and CGSG-SSR which solved nearly all trials. At $\overline{\epsilon}=1.0$, none of the compared methods solved any trials on P5. On P3, SAGE fared comparably to CFD-SSR, whilst FFD-SSR achieved a substantially better data profile at $\overline{\epsilon}=1\times10^{-3}$. High noise on P3 ($\overline{\epsilon}=1.0$) proved very difficult even for the noise-informed estimators, with all compared methods achieving less than 20\% solved trials by the evaluation budget.

Taken together, SAGE's noise-adaptive sampling (without \textit{a priori} knowledge of $\overline{\epsilon}$) did not uniformly outperform noise-informed competitors: it is markedly less efficient than FFD-SSR on P1-P2, comparable to most competitors on P4 and P5, and outperformed on P3. It remains, however, the only estimator that achieved any of these results without $\overline{\epsilon}$ being supplied \textit{a priori}, and was never the worst performer across the five problems.

However, in practical settings with expensive and time-consuming evaluations (e.g., computational fluid dynamics simulations, or automotive experiments) the noise characterisation to obtain the noise-informed step sizes for FFD-SSR, CFD-SSR, GSG-SSR, and CGSG-SSR is itself often infeasible. Noise characterisation requires an initial data set dedicated to estimating $\overline{\epsilon}$, which is exactly the kind of expensive prerequisite these applications cannot easily afford. Therefore, the current comparison is only a sanity check against \textit{noise-informed} gradient estimators. Even without this \textit{a priori} information, SAGE reaches a comparable final success rate to the noise-informed estimators on P1, P2, and P4 within the evaluation budget (albeit less efficiently than FFD-SSR on P1-P2), and is never the single worst performer on any of the five problems. This combination of adaptive sampling and samples reuse, without upfront noise characterisation, is therefore particularly promising for the expensive-evaluation settings that motivate this work, where Figs.~\ref{fig:data-profiles-lsq-varying-noise}-\ref{fig:data-profiles-l2logreg-varying-noise} (no noise information for any method) reflect the more realistic operating condition, and where SAGE's advantage is most pronounced.

\section{Conclusion}
\label{sec:conclusion}

In this paper, we pose gradient estimation from zeroth-order information as a set membership problem rather than a point-estimation one, thereby reversing the conventional estimator-first paradigm in which a single gradient estimate is computed, and its accuracy analysed only post-hoc. Instead, we first characterise the set of admissible gradients consistent with the available samples and regularity assumptions, making uncertainty directly usable for deciding whether to refine the gradient information, or extract a representative point gradient. This enables the natural (re)use of global information contained in past samples, and avoids reliance on fixed sampling heuristics. Within this framework, our proposed Set-Based Adaptive Gradient Estimator (SAGE) turns this admissible set into a practical, sample-efficient algorithm: it bootstraps an initial estimate, refines it only when the set's geometry indicates that further samples are informative, and extracts a deterministic point estimate via a rotated bounding hyperbox: all without requiring \textit{a priori} knowledge of the function's smoothness or noise characteristics. Our results naturally extend to noisy evaluations with unknown finite bounds, where we show that with noise, the achievable gradient precision is fundamentally limited, and best described by a set rather than a single estimate.

Our benchmark results confirm that this comes with a practically relevant trade-off. In the gradient accuracy benchmark, SAGE is competitive with, or better than, all compared estimators once noise is present, despite operating without \textit{a priori} noise information. In the first-order optimisation benchmark, this advantage is more nuanced. In particular, when the competing finite-difference and smoothing-based estimators are given \textit{a priori} noise characterisation and their step sizes adjusted accordingly, SAGE is less efficient than the best-informed competitor on some problems, but in general broadly comparable with other noise-informed methods. Nevertheless, SAGE remains the only estimator among those compared that achieves these results without upfront noise characterisation. This is precisely the setting we target: in industrial applications such as computational fluid dynamics simulations or automotive experiments, function evaluations are typically expensive and time-consuming, making the initial data collection needed for noise characterisation itself impractical. In such settings, SAGE's combination of adaptive sampling and sample reuse offers a favourable trade-off between accuracy and the practical cost of deployment, without requiring information that is often unavailable in practice.

\section*{Data Availability}
The authors declare that the source code used to generate the data analysed in the comparative study in this paper is available online at \url{https://github.com/lorenzosabugjr/SAGE}. Other raw data files, if required in another format, are available from the corresponding author upon reasonable request.

\begin{acknowledgements}
This research has been supported by the Italian Ministry of Enterprises and Made in Italy in the framework of the project 4DDS (4D Drone Swarms) under grant no. F/310097/01-04/X56, and by the Italian Ministry of University and Research under grants ``DeepAirborne--Advanced Modeling, Control and Design Optimization Methods for Deep Offshore Airborne Wind Energy'' (NextGenerationEU fund, project P2022927H7), and ``Learning-based Model Predictive Control by Exploration and Exploitation in Uncertain Environments'' (PRIN PNRR 2022 fund, ID P2022EXP2W)
\end{acknowledgements}

\appendix  

\section*{Appendix: Proofs}

\begin{proof}[Corollary~\ref{corollary:gradient-accuracy-2samples}]
We directly invoke Lemma 4.1.1 of Nesterov \cite{nesterov2018lectures}, which states that

\begin{equation}
\label{eqn:nesterov2018-equation}
     \left|\idxval{z}{j}{} - \idxval{z}{i}{} - \idxptT{g}{i}{}\bm{a}_{ij} - \frac{1}{2}{} \bm{a}_{ij}^\top \idxpt{H}{i}{} \bm{a}_{ij} \right| \leq \frac{1}{6} \gamma_H \|\bm{a}_{ij}\|^3.
\end{equation}

Dividing both sides by $\|\bm{a}_{ij}\|$,
\[
     \left|\frac{\idxval{z}{j}{} - \idxval{z}{i}{}}{\|\bm{a}_{ij}\|} - \idxptT{g}{i}{}\bm{u}_{ij} - \frac{1}{2}{} \|\bm{a}_{ij}\|\bm{u}_{ij}^\top \idxpt{H}{i}{} \bm{u}_{ij} \right| \leq \frac{1}{6} \gamma_H \|\bm{a}_{ij}\|^2.
\]

\noindent Applying the inequality $|\bm{u}_{ij}^\top \idxpt{H}{i}{} \bm{u}_{ij}| \leq \idxval{H}{i}{}$, we now have

\[
     \Big|\frac{\idxval{z}{j}{} - \idxval{z}{i}{}}{\|\bm{a}_{ij}\|} - \idxptT{g}{i}{}\bm{u}_{ij} \Big| \leq  \frac{1}{2}{} \idxval{H}{i}{} \|\bm{a}_{ij}\| + \frac{1}{6}\gamma_H \|\bm{a}_{ij}\|^2.
\]

\noindent which results in \eqref{eqn:gradient-slab}, concluding the proof.
\qed
\end{proof}

\begin{proof}[Theorem~\ref{theorem:gradient-set}]

Equation~\eqref{eqn:gradient-slab} becomes

\[
    -\frac{1}{2} \|\bm{a}_{ij}\| \idxval{H}{i}{} - \frac{1}{6} \|\bm{a}_{ij}\|^2 \gamma_H \leq \tilde{s}_{ij} - \bm{g}^{(i)\top} \bm{u}_{ij} \leq \frac{1}{2} \|\bm{a}_{ij}\| \idxval{H}{i}{} + \frac{1}{6} \|\bm{a}_{ij}\|^2 \gamma_H.
\]

Splitting the above, we now build the following inequalities

\begin{align}
    \tilde{s}_{ij} - \bm{g}^{(i)\top} \bm{u}_{ij} & \leq \frac{1}{2} \|\bm{a}_{ij}\| \idxval{H}{i}{} + \frac{1}{6} \|\bm{a}_{ij}\|^2 \gamma_H \\
    -\tilde{s}_{ij} + \bm{g}^{(i)\top} \bm{u}_{ij} & \leq \frac{1}{2} \|\bm{a}_{ij}\| \idxval{H}{i}{} + \frac{1}{6} \|\bm{a}_{ij}\|^2 \gamma_H
\end{align}

\noindent and then

\begin{align}
     -\bm{u}_{ij}^\top \bm{g}^{(i)} - \frac{1}{2}\|\bm{a}_{ij}\|\idxval{H}{i}{} - \frac{1}{6} \|\bm{a}_{ij}\|^2 \gamma_H & \leq -\tilde{s}_{ij} \\
      \bm{u}_{ij}^\top \bm{g}^{(i)} - \frac{1}{2}\|\bm{a}_{ij}\|\idxval{H}{i}{} - \frac{1}{6} \|\bm{a}_{ij}\|^2 \gamma_H & \leq  \tilde{s}_{ij}    
\end{align}

Reorganising to matrix inequalities, we get

\begin{equation}
    \begin{bmatrix}
        -\bm{u}_{ij}^\top & -\frac{1}{2}\|\bm{a}_{ij}\| & -\frac{1}{6}\|\bm{a}_{ij}\|^2 \\
         \bm{u}_{ij}^\top & -\frac{1}{2}\|\bm{a}_{ij}\| & -\frac{1}{6}\|\bm{a}_{ij}\|^2 \\
    \end{bmatrix}
    \begin{bmatrix}
        \bm{g}^{(i)} \\ \idxval{H}{i}{} \\ \gamma_H
    \end{bmatrix}
    \leq
    \begin{bmatrix}
        -\tilde{s}_{ij} \\ \tilde{s}_{ij}
    \end{bmatrix} \label{eqn:gradient-est-ineq}
\end{equation}

Now, collecting such matrix inequalities for all $\idxpt{x}{j}{} \in \itrpt{X}{n}{} \setminus \idxpt{x}{i}{}$, we get \eqref{eqn:gradient-set}, and the theorem is proven.
\qed
\end{proof}

\begin{proof}[Corollary~\ref{corollary:gradient-accuracy-2samples-noise}]
    Given two samples $(\idxpt{x}{i}{},\idxval{z}{i}{})$ and $(\idxpt{x}{j}{},\idxval{z}{j}{})$ and recognising that $f(\idxpt{x}{i}{}) = \idxval{z}{i}{} - \epsilon_i$, Lemma 4.1.1 of \cite{nesterov2018lectures} given as \eqref{eqn:nesterov2018-equation} becomes

    \[
         \left|\idxval{z}{j}{} - \idxval{z}{i}{} + \epsilon_i - \epsilon_j - \idxptT{g}{i}{}\bm{a}_{ij} - \frac{1}{2}{} \bm{a}_{ij}^\top \idxpt{H}{i}{} \bm{a}_{ij} \right| \leq \frac{1}{6} \gamma_H \|\bm{a}_{ij}\|^3.
    \]

    \noindent Following the same steps as in the proof of Corollary~\ref{corollary:gradient-accuracy-2samples} and recognising that $\frac{\epsilon_j - \epsilon_i}{\|\bm{a}_{ij}\|} \leq \frac{2\overline{\epsilon}}{\|\bm{a}_{ij}\|}$, we arrive at \eqref{eqn:directional-derivative-uncertainty-noise}, which completes the proof.
    \qed
\end{proof}


\end{document}